\documentclass[10pt]{article}
\usepackage{amsmath, amsthm, amssymb, verbatim, tikz, xcolor, subfigure}
\usepackage{latexsym}
\usepackage{amsfonts}
\usepackage{graphicx}

\usepackage{float}
\usepackage{bigints}
\usepackage{fullpage}

\usepackage{footmisc}
\DeclareMathOperator{\coin}{coin}
\DeclareMathOperator{\ill}{ill}
\DeclareMathOperator{\conv}{conv}
\DeclareMathOperator{\vol}{vol}
\DeclareMathOperator{\wcoin}{coin_{w}}
\DeclareMathOperator{\cyl}{cyl}

\newcommand{\crv}{{\rm crv}}
\newcommand{\Eu}{\mathbb{E}}

\newcommand{\KK}{\mathbf{K}}

\begin{document}

\newtheorem{theorem}{Theorem}[section]
\newtheorem{question}[theorem]{Question}
\newtheorem{example}[theorem]{Example}
\newtheorem{observation}[theorem]{Observation}
\newtheorem{remark}[theorem]{Remark}
\newtheorem{fact}{Fact}
\newtheorem{proc}[theorem]{Procedure}
\newtheorem{conjecture}[theorem]{Conjecture}
\newtheorem{lemma}[theorem]{Lemma}
\newtheorem{proposition}[theorem]{Proposition}
\newtheorem{corollary}[theorem]{Corollary}
\newtheorem{fac}[theorem]{Fact}
\newtheorem*{abst}{Abstract}
\newtheorem{definition}{Definition}
\newtheorem{problem}{Problem}
\newtheorem{algo}{Algorithm}

\title{The geometry of homothetic covering and illumination} 

\author{K\'{a}roly Bezdek and Muhammad A. Khan}
\date{}
 \maketitle

\begin{abstract}
At a first glance, the problem of illuminating the boundary of a convex body by external light sources and the problem of covering a convex body by its smaller positive homothetic copies appear to be quite different. They are in fact two sides of the same coin and give rise to one of the important longstanding open problems in discrete geometry, namely, the Illumination Conjecture. In this paper, we survey the activity in the areas of discrete geometry, computational geometry and geometric analysis motivated by this conjecture. Special care is taken to include the recent advances that are not covered by the existing surveys. We also include some of our recent results related to these problems and describe two new approaches -- one conventional and the other computer-assisted -- to make progress on the illumination problem. Some open problems and conjectures are also presented.            
\vspace{2mm}

\noindent \textit{Keywords and phrases:}  Illumination number, Illumination Conjecture, Covering Conjecture, Separation Conjecture, X-ray number, X-ray Conjecture, illumination parameter, covering parameter, covering index, cylindrical covering parameters, $\epsilon$-net of convex bodies. 

\vspace{2mm}

\noindent \textit{MSC (2010):} 52A37, 52A40, 52C15, 52C17. 
\end{abstract}

\section{Shedding some `light'}\label{sec:intro}
\begin{quotation}
\small \textit{``$\ldots$ $N_{k}$  bezeichne die kleinste nat\"urliche Zahl, f\"ur welche die nachfolgende Aussage richtig ist: Ist $A$ ein eigentlicher konvexer K\"{o}rper des $k$-dimensionalen euklidischen Raumes, so gibt es $n$ mit $A$ translations-gleiche K\"{o}rper $A_i$ mit $n \leq N_{k}$ derart, dass jeder Punkt von $A$ ein innerer Punkt der Vereinigungsmenge $\bigcup_{i}A_{i}$ ist,$\ldots$Welchen Wert hat $N_{k}$ f\"{u}r $k\geq 3$?''} \cite{hadwiger1}
\end{quotation}
\normalsize
The above statement roughly translates to ``Let $N_{k}$ denote the smallest natural number such that any $k$-dimensional convex body can be covered by the interior of a union of at the most $N_{k}$ of its translates. What is $N_{k}$ for $k\geq 3$?'' When Hadwiger raised this question in 1957 he probably did not imagine that it would remain unresolved half a century later and become a central problem in discrete geometry. Apparently, Hadwiger had a knack of coming up with such questions\footnote{The Hadwiger conjecture in graph theory is, in the words of Bollob\'{a}s, Catalin and Erd\H{o}s \cite{bollobas1}, ``\textit{one of the deepest unsolved problems in graph theory}". Hadwiger even edited a column on unsolved problems in the journal \textit{Elemente der Mathematik}. On the occasion of Hadwiger's 60th birthday, Victor Klee dedicated the first article in the Research Problems section of the {\it American Mathematical Monthly} to Hadwiger's work on promoting research problems \cite[pp. 389--390]{gardner}.}. However, he was not the first one to study this particular problem. In fact, its earliest occurrence can be traced back to Levi's 1955 paper \cite{levi1}, who formulated and settled the 2-dimensional case of the problem. Later in 1960, the question was restated by Gohberg and Markus\footnote{Apparently, Gohberg and Markus worked on the problem independently without knowing about the work of Levi and Hadwiger \cite{boltyanski-gohberg}.} \cite{gohberg1} in terms of covering by homothetic copies. The equivalence of both formulations is relatively easy to check and details appear in Section 34 of \cite{boltyanski2}. 

\begin{conjecture}[\textbf{Covering Conjecture}]\label{hadwigerconjecture}
We can cover any $d$-dimensional convex body by $2^d$ or fewer of its smaller positive homothetic copies in Euclidean $d$-space, $d\ge 3$. Furthermore, $2^{d}$ homothetic copies are required only if the body is an affine $d$-cube. 
\end{conjecture} 

The same conjecture has also been referred to in the literature as the Levi--Hadwiger Conjecture, Gohberg--Markus Covering Conjecture and Hadwiger Covering Conjecture. The condition $d\geq 3$ has been added as the statement is known to be true in the plane \cite{levi1, hadwiger2}.

Let us make things formal. A $d$-dimensional \textit{convex body} $\mathbf{K}$ is a compact convex subset of the Euclidean $d$-space, ${\mathbb{E}}^{d}$ with nonempty interior. Let $\mathbf{o}$ denote the origin of ${\mathbb{E}}^{d}$. Then $\mathbf{K}$ is said to be \textit{$\mathbf{o}$-symmetric} if $\mathbf{K}=-\mathbf{K}$ and \textit{centrally symmetric} if some translate of $\mathbf{K}$ is $\mathbf{o}$-symmetric. Since the quantities studied in this paper are invariant under affine transformations, we use the terms $\mathbf{o}$-symmetric and centrally symmetric interchangeably. A \textit{homothety} is an affine transformation of ${\mathbb{E}}^{d}$ of the form $\mathbf{x}\mapsto \mathbf{t}+\lambda\mathbf{x}$, where $\mathbf{t}\in {\mathbb{E}}^{d}$ and $\lambda$ is a non-zero real number. The image $\mathbf{t}+\lambda \mathbf{K}$ of a convex body $\mathbf{K}$ under a homothety is said to be its \textit{homothetic copy} (or simply a \textit{homothet}). A homothetic copy is \textit{positive} if $\lambda >0$ and \textit{negative} otherwise. Furthermore, a homothetic copy with $0<\lambda<1$ is called a smaller positive homothet. In terms of the notations just introduced the Covering Conjecture states that for any $\mathbf{K}\subseteq {\mathbb{E}}^{d}$, there exist $\mathbf{t}_{i}\in {\mathbb{E}}^{d}$ and $0<\lambda_i<1$, for $i=1,\ldots,2^{d}$, such that 
\begin{equation}\label{maineq}
\mathbf{K}\subseteq \bigcup_{i=1}^{2^{d}}(\mathbf{t}_{i}+\lambda_{i}\mathbf{K}). 
\end{equation} 

\vspace{-6mm}
\begin{figure}[htb]
	\centering
		\includegraphics[scale=0.4]{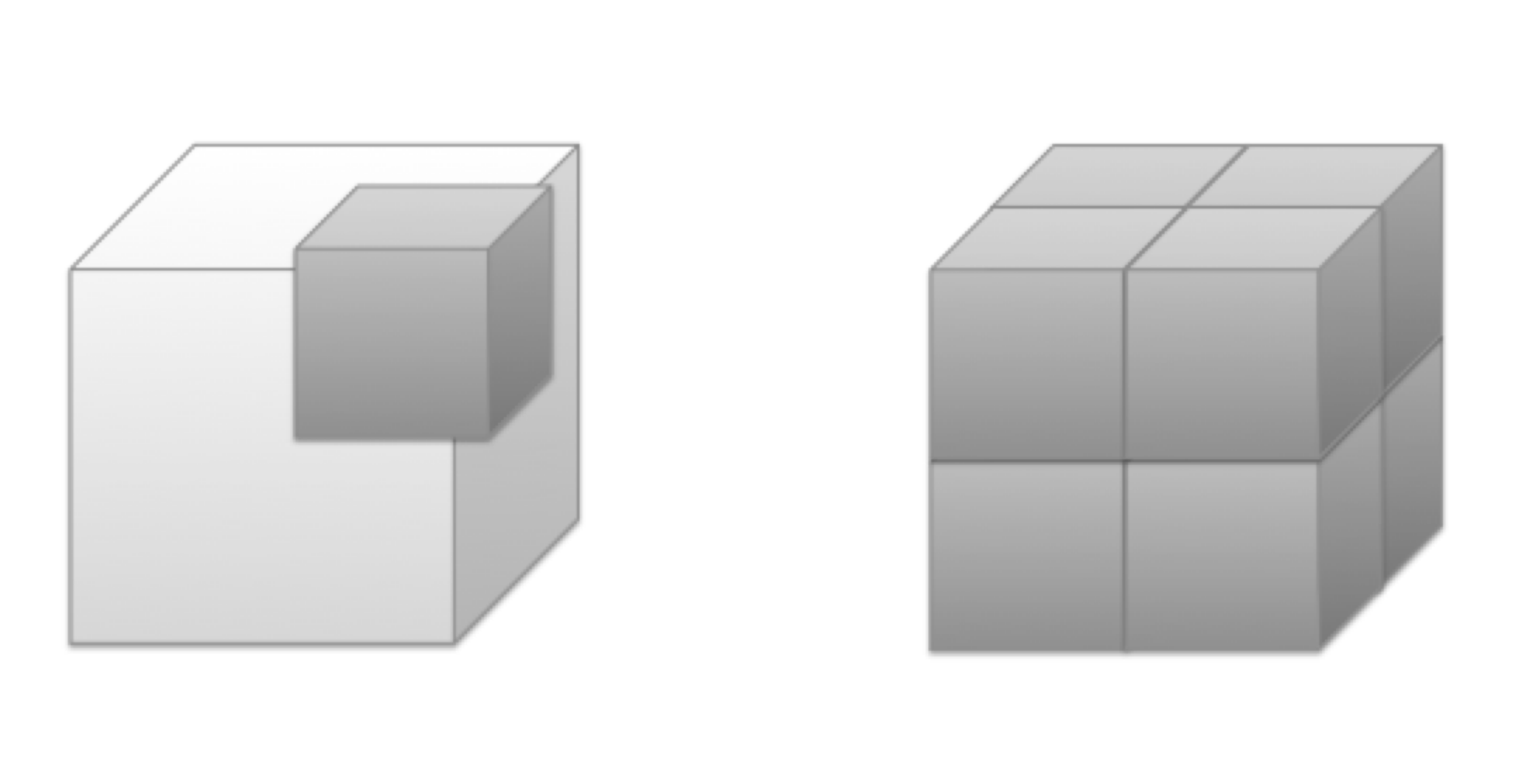}
	\label{fig:homothetic}
\vspace{-6mm}
	\caption{A cube can be covered by 8 smaller positive homothets and no fewer.}
\end{figure}

A light source at a point $\mathbf{p}$ outside a convex body $\mathbf{K}\subset{\mathbb{E}}^{d}$, \textit{illuminates} a point $\mathbf{x}$ on the boundary of $\mathbf{K}$ if the halfline originating from $\mathbf{p}$ and passing through $\mathbf{x}$ intersects the interior of $\mathbf{K}$ at a point not lying between $\mathbf{p}$ and $\mathbf{x}$. The set of points $\{\mathbf{p}_{i}:i=1,\ldots,n\}$ in the exterior of $\mathbf{K}$ is said to \textit{illuminate} $\mathbf{K}$ if every boundary point of $\mathbf{K}$ is illuminated by some $\mathbf{p}_{i}$. The \textit{illumination number} $I(\mathbf{K})$ of $\mathbf{K}$ is the smallest $n$ for which $\mathbf{K}$ can be illuminated by $n$ point light sources.

\begin{figure}[htb]
	\label{fig:points}
	\centering
\begin{tabular}{ccc}
\subfigure [ ]{
\begin{tikzpicture}[scale = 0.6, baseline]
\large
\draw[fill=black!20, rotate=30, very thick] (0,0) ellipse [x radius=2cm, y radius=1cm];
\draw[fill=black] (-1.1,2.5) circle [radius=2pt] node[above] {$\mathbf{p}$};
\draw[fill=black] (-0.1,1.1) circle [radius=2pt] node[above] {$\mathbf{x}$};
\draw[->, very thick] (-1.1,2.5) -- (0.5,0.15);
\normalsize
\end{tikzpicture}
} &    &
\subfigure [ ]{
\begin{tikzpicture}[scale=0.6, baseline]
\large
\draw[fill=black!20, rotate=30, very thick] (0,0) ellipse [x radius=2cm, y radius=1cm];
\draw[fill=black] (-1.1,2.5) circle [radius=2pt] node[above] {$\mathbf{p_1}$};
\draw[->, very thick] (-1.1,2.5) -- (1.3,1.3); 
\draw[->, very thick] (-1.1,2.5) -- (-0.2,1); 
\draw[->, very thick] (-1.1,2.5) -- (-1.1,0.5); 
\draw[->, very thick] (-1.1,2.5) -- (-1.8,-0.55); 

\draw[fill=black] (4,0.3) circle [radius=2pt] node[above] {$\mathbf{p_2}$};
\draw[->, very thick] (4,0.3) -- (1.3,1.3); 
\draw[->, very thick] (4,0.3) -- (1.85,.55); 
\draw[->, very thick] (4,0.3) -- (0.7,-.8); 

\draw[fill=black] (-2.8,-3) circle [radius=2pt] node[left] {$\mathbf{p_3}$};
\draw[->, very thick] (-2.8,-3) -- (-1.8,-0.55); 
\draw[->, very thick] (-2.8,-3) -- (-1.5,-1.2); 
\draw[->, very thick] (-2.8,-3) -- (-.5,-1.3); 
\draw[->, very thick] (-2.8,-3) -- (0.7,-.8); 
\normalsize
\end{tikzpicture}
}
\end{tabular}
\caption{(a) Illuminating a boundary point $\mathbf{x}$ of $\mathbf{K}\subset{\mathbb{E}}^{d}$ by the point light source $\mathbf{p}\in{\mathbb{E}}^{d}\setminus\mathbf{K} $, (b) $I(K) = 3$.}
\end{figure}
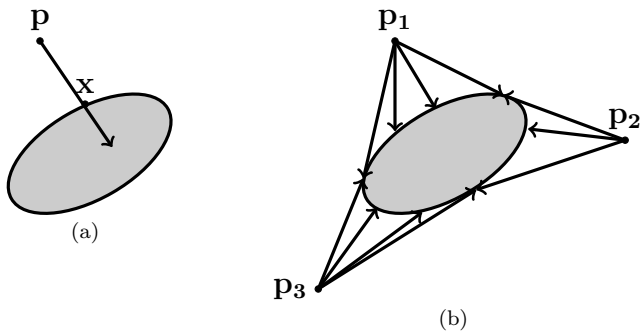

One can also consider illumination of $\mathbf{K}\subset{\mathbb{E}}^{d}$ by parallel beams of light. Let $\mathbb{S}^{d-1}$ be the unit sphere centered at the origin $\mathbf{o}$ of $\mathbb{E}^{d}$. We say that a point $\mathbf{x}$ on the boundary of $\mathbf{K}$ is illuminated in the direction $\mathbf{v}\in\mathbb{S}^{d-1} $ if the halfline originating from $\mathbf{x}$ and with direction vector $\mathbf{v}$ intersects the interior of $\mathbf{K}$. 

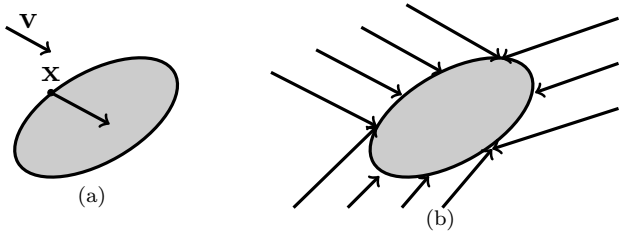
\begin{figure}[htb]
	\label{fig:lines}
	\centering
\begin{tabular}{ccc}
\subfigure [ ]{
\begin{tikzpicture}[scale = 0.6, baseline]
\large
\draw[fill=black!20, rotate=30, very thick] (0,0) ellipse [x radius=2cm, y radius=1cm];
\draw[fill=black] (-1,0.55) circle [radius=2pt];
\draw[->, very thick] (-1,0.55) node[above] {$\mathbf{x}$} -- (0.3,-.15); 
%\draw[->, very thick] (-2,2.05) -- (-0.7,1.5) node[midway,above] {$\vec{\ell}$}; 
%\draw[->, very thick] (-2,2) -- (-0.2,1);  
\draw[->, very thick] (-2,2) -- (-1,1.45) node[midway,above] {$\mathbf{v}$};  
\normalsize
\end{tikzpicture}
} &    &
\subfigure [ ]{
\begin{tikzpicture}[scale = 0.6, baseline]
\large
\draw[fill=black!20, rotate=30, very thick] (0,0) ellipse [x radius=2cm, y radius=1cm];
\draw[->, very thick] (-1,2.5) -- (1.1,1.3); 
\draw[->, very thick] (-2,2) -- (-0.2,1); 
\draw[->, very thick] (-3,1.5) -- (-1.1,0.5); 
\draw[->, very thick] (-4,1) -- (-1.65,-0.2); 

\draw[->, very thick] (3.7,2.2) -- (1.1,1.3); 
\draw[->, very thick] (3.7,1.2) -- (1.85,.55); 
\draw[->, very thick] (3.7,.2) -- (0.9,-.7); 

\draw[->, very thick] (-3.5,-2) -- (-1.65,-0.2); 
\draw[->, very thick] (-2.3,-2) -- (-1.6,-1.3); 
\draw[->, very thick] (-1.1,-2) -- (-.5,-1.3); 
\draw[->, very thick] (-.2,-2) -- (0.9,-.7); 
\normalsize
\end{tikzpicture}
}
\end{tabular}
\caption{(a) Illuminating a boundary point $\mathbf{x}$ of $\mathbf{K}\subset{\mathbb{E}}^{d}$ by a direction $\mathbf{v}\in\mathbb{S}^{d-1} $, (b) $I(K) = 3$.}
\end{figure}

The former notion of illumination was introduced by Hadwiger \cite{hadwiger2}, while the latter notion is due to Boltyanski\footnote{Vladimir Boltyanski (also written Boltyansky, Boltyanskii and Boltjansky) is a prolific mathematician and recepient of Lenin Prize in science. He has authored more than 220 mathematical works including, remarkably, more than 50 books!}\cite{boltyanski1}. It may not come as a surprise that the two concepts are equivalent in the sense that a convex body $\mathbf{K}$ can be illuminated by $n$ point sources if and only if it can be illuminated by $n$ directions. However, it is less obvious that any covering of $\mathbf{K}$ by $n$ smaller positive homothetic copies corresponds to illuminating $\mathbf{K}$ by $n$ points (or directions) and vice versa (see \cite{boltyanski2} for details). Therefore, the following Illumination Conjecture \cite{boltyanski1, boltyanski2, hadwiger2} of Hadwiger and Boltyanski is equivalent to the Covering Conjecture. 

\begin{conjecture}[\textbf{Illumination Conjecture}] 
The illumination number $I(\mathbf{K})$ of any $d$-dimensional convex body $\mathbf{K}$, $d\geq 3$, is at most $2^d$ and $I(K) = 2^d$ only if $\mathbf{K}$ is an affine $d$-cube.
\end{conjecture}

\vspace{-4mm}
\begin{figure}[!h]
\centering 
\includegraphics[scale=0.5]{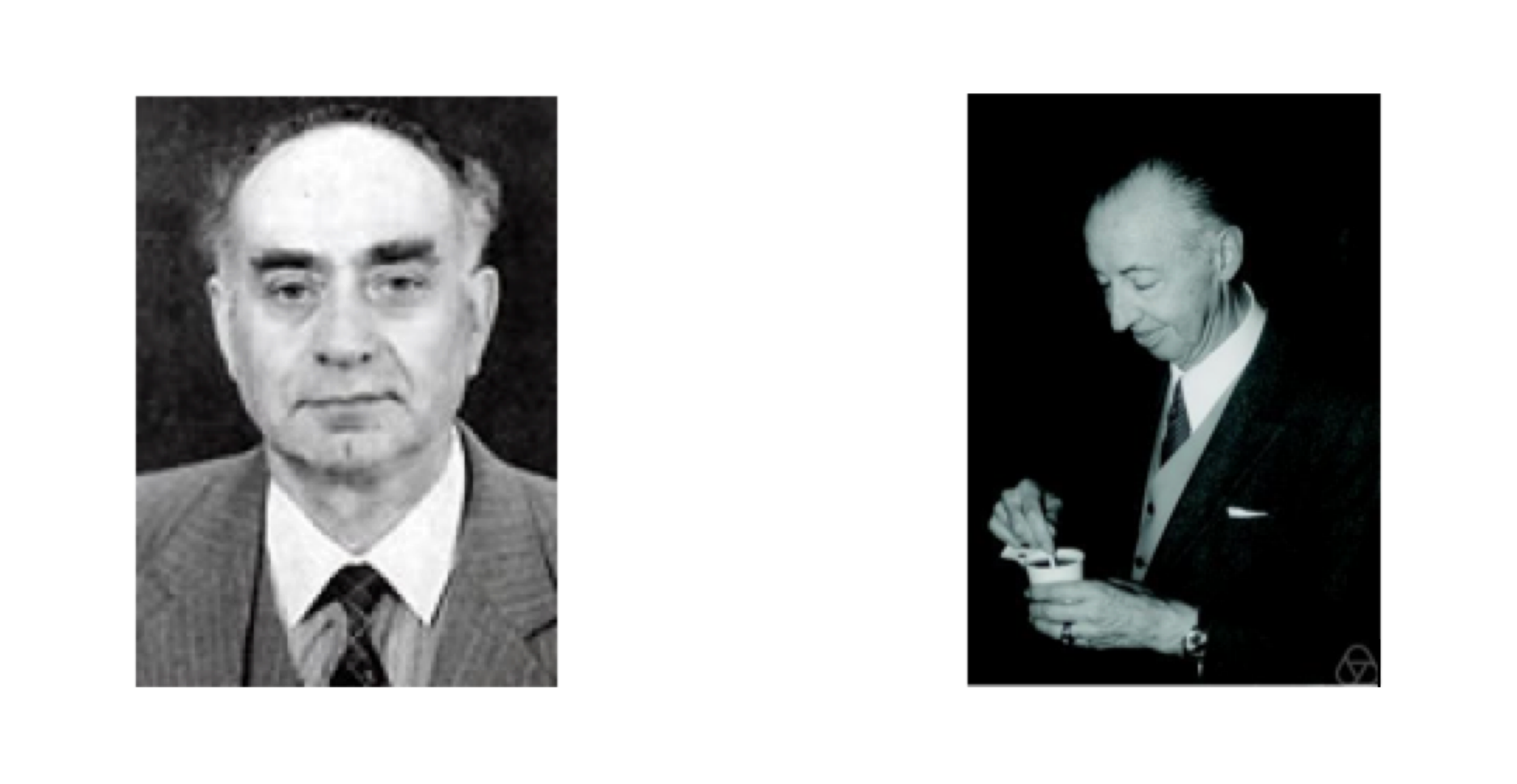}
\vspace{-4mm} 
\caption{Vladimir Boltyanski (left, courtesy Annals of the Moscow University) and Hugo Hadwiger (right, courtesy Oberwolfach Photo Collection), two of the main proponents of the illumination problem.}
\end{figure}

The conjecture also asserts that affine images of $d$-cubes are the only extremal bodies. The conjectured bound of $2^{d}$ results from the $2^{d}$ vertices of an affine cube, each requiring a different light source to be illuminated. In the sequel, we use the titles Covering Conjecture and Illumination Conjecture interchangeably, shifting between the covering and illumination paradigms as convenient.    

We have so far seen three equivalent formulations of the Illumination Conjecture. But there are more. In fact, it is perhaps an indication of the richness of this problem that renders it to be studied from many angles, each with its own intuitive significance. We state one more equivalent form found independently by P. Soltan and V. Soltan \cite{soltan1}, who formulated it for the $\mathbf{o}$-symmetric case only and the first author \cite{bezdek-conj1, bezdek-conj2}. 

\begin{conjecture}[\textbf{Separation Conjecture}]
Let $\mathbf{K}$ be an arbitrary convex body in ${\mathbb{E}}^d$, $d \geq 3$, and $\mathbf{o}$ be an arbitrary interior point of $\mathbf{K}$. Then there exist $2^d$ hyperplanes of ${\mathbb{E}}^d$ such
that each intersection of $\mathbf{K}$ with a supporting hyperplane, called a \textit{face} of $\mathbf{K}$, can be strictly separated from $\mathbf{o}$ by at least one of the $2^d$ hyperplanes. Furthermore, $2^d$ hyperplanes are needed only if $\mathbf{K}$ is the convex hull of $d$ linearly independent line segments which intersect at the common relative interior point $\mathbf{o}$.
\end{conjecture}  

Over the years, the illumination conjecture has inspired a vast body of research in convex and discrete geometry, computational geometry and geometric analysis. There exist some nice surveys on the topic such as the papers \cite{bezdek-survey, martini1} and the corresponding chapters of the books \cite{bezdek-book, boltyanski2}. However, most of these are a bit dated. Moreover, we feel that the last few years have seen some interesting new ideas, such as the possibility of a computer-assisted proof, that are not covered by any of the above-mentioned surveys. The aim of this paper is to provide an accessible introduction to the geometry surrounding the Illumination Conjecture and a snapshot of the research motivated by it, with special emphasis on some of the recent developments. At the same time we describe some of our new results in this area.   

We organize the material as follows. Section \ref{illumination-progress} gives a brief overview on the progress of the Illumination Conjecture. In Section \ref{sec:relatives}, we mention some important relatives of the illumination problem, while Section \ref{sec:quantify} explores the known important quantitative versions of the problem including a new approach to make progress on the Illumination Conjecture based on the covering index of convex bodies (see Problem \ref{ball-max} and the discussion following it in Section \ref{sec:index}). Finally, in Section \ref{sec:computational} we present Zong's computer-assisted approach \cite{zong1} for possibly resolving the Illumination Conjecture in low dimensions.

\section{Progress on the Illumination Conjecture}\label{illumination-progress}
\subsection{Results in ${\mathbb{E}}^{3}$ and ${\mathbb{E}}^{4}$} \label{sec:E3-E4}
%\subsection*{Results in dimension 3}

Despite its intuitive richness, the illumination conjecture has been notoriously difficult to crack even in the first nontrivial case of $d=3$. The closest anything has come is the proof announced by Boltyanski \cite{boltyanski3} for the 3-dimensional case. Unfortunately, the proof turned out to have gaps that remain to date. Later, Boltyanski \cite{boltyanski6} modified his claim to the following.  

\begin{theorem}\label{boltyanski-modified}
Let $\mathbf{K}$ be a convex body of ${\mathbb{E}}^{3}$ with ${\rm md \ }\mathbf{K} = 2$. Then $I(\mathbf{K})\leq 6$. 
\end{theorem}

Here ${\rm md}$ is a functional introduced by Boltyanski in \cite{boltyanski-md} and defined as follows for any $d$-dimensional convex body: Let $\mathbf{K} \subseteq {\mathbb{E}}^{d}$ be a convex body. Then ${\rm md}(\mathbf{K})$ is the greatest integer $m$ for which there exist $m + 1$ regular boundary points of $\mathbf{K}$ such that the outward unit normals $\mathbf{v}_{0},\ldots , \mathbf{v}_{m}$ of $\mathbf{K}$ at these points are minimally dependent, i.e., they are the vertices of an $m$-dimensional simplex that contains the origin in its relative interior.\footnote{In fact, it is proved in \cite{boltyanski-md} that ${\rm md}(\mathbf{K})={\textnormal{him}}(\mathbf{K})$ holds for any convex body $\mathbf{K}$ of ${\mathbb{E}}^{d}$ and therefore one can regard Theorem~\ref{boltyanski-modified} as an immediate corollary of Theorem~\ref{him} for $d=3$ in Section~\ref{sec:general}.}

So far the best upper bound on illumination number in three dimensions is due to Papadoperakis \cite{papadoperakis1}. 

\begin{theorem}
The illumination number of any convex body in ${\mathbb{E}}^3$ is at most 16. 
\end{theorem} 

However, there are partial results that establish the validity of the conjecture for some large classes of convex bodies. Often these classes of convex bodies have some underlying symmetry. Here we list some such results. A convex polyhedron $\mathbf{P}$ is said to have affine symmetry if the affine symmetry group of $\mathbf{P}$ consists of the identity and at least one other affinity of ${\mathbb{E}}^3$. The first author obtained the following result \cite{bezdek-conj1}. 

\begin{theorem}\label{bezdektheorem}
If $\mathbf{P}$ is a convex polyhedron of ${\mathbb{E}}^3$ with affine symmetry, then the illumination number of $\mathbf{P}$ is at most 8. 
\end{theorem}

Recall that a convex body $\mathbf{K}$ is said to be centrally symmetric if it has a point of symmetry. Furthermore, a body $\mathbf{K}$ is symmetric about a plane $p$ if a reflection across that plane leaves $\mathbf{K}$ unchanged. Lassak \cite{lassak1} proved that under the assumption of central symmetry, the illumination conjecture holds in three dimensions. 

\begin{theorem}\label{Lassak-illumination}
If $\mathbf{K}$ is a centrally symmetric convex body in ${\mathbb{E}}^3$, then $I(\mathbf{K})\leq 8$.
\end{theorem}   

Dekster \cite{dekster1} extended Theorem \ref{bezdektheorem} from polyhedra to convex bodies with plane symmetry. 

\begin{theorem}
If $\mathbf{K}$ is a convex body symmetric about a plane in ${\mathbb{E}}^3$, then $I(\mathbf{K})\leq 8$.
\end{theorem}

It turns out that for 3-dimensional bodies of constant width -- that is bodies whose width, measured by the distance between two opposite parallel hyperplanes touching its boundary, is the same regardless of the direction of those two parallel planes --  we get an even better bound. 

\begin{theorem}
The illumination number of any convex body of constant width in ${\mathbb{E}}^3$ is at most 6. 
\end{theorem} 

Proofs of the above theorem have appeared in several papers \cite{bezdek-langi1, lassak2,weissbach1}. It is, in fact, reasonable to conjecture the following even stronger result. 

\begin{conjecture}
The illumination number of any convex body of constant width in ${\mathbb{E}}^3$ is exactly 4. 
\end{conjecture}

The above conjecture, if true, would provide a new proof of Borsuk's conjecture \cite{borsuk} in dimension three, which states that any set of unit diameter in ${\mathbb{E}}^{3}$ can be partitioned into at most four subsets of diameter less than one.  We remark that although it is false in general \cite{kahn}, Borsuk's conjecture has a long and interesting history of its own and the reader can look up \cite{boltyanski4, boltyanski2, grunbaum} for detailed discussions. 

Now let us consider the state of the Illumination Conjecture in ${\mathbb{E}}^4$. It is well known that neighbourly $d$-polytopes have the maximum number of facets among $d$-polytopes with a fixed number of vertices (for more details on this see for example, \cite{BiFo}). Thus, it is natural to investigate the Separation Conjecture for neighbourly $d$-polytopes (see also Theorem~\ref{dual}). Since interesting neighbourly $d$-polytopes exist only in $\Eu^d$ for $d\geq4$, it is particularly natural to first restrict our attention to neighbourly $4$-polytopes. Starting from a cyclic $4$-polytope, the sewing procedure of Shemer (for details see \cite{BiFo}) produces an infinite family of neighbourly $4$-polytopes each of which is obtained from the previous one by adding one new vertex in a suitable way. Neighbourly $4$-polytopes obtained from a cyclic $4$-polytope by a sequence of sewings are called {\it totally-sewn}. The main result of the very recent paper \cite{BiFo} of Bisztriczky and Fodor is a proof of the Separation Conjecture for totally-sewn neighbourly $4$-polytopes. 

\begin{theorem}
Let $\mathbf{P}$ be an arbitrary totally-sewn neighbourly $4$-polytope in ${\mathbb{E}}^4$, and $\mathbf{o}$ be an arbitrary interior point of $\mathbf{P}$. Then there exist $16$ hyperplanes of ${\mathbb{E}}^4$ such that each face of $\mathbf{P}$, can be strictly separated from $\mathbf{o}$ by at least one of the $16$ hyperplanes.
\end{theorem}

However, Bisztriczky \cite{Bi02} conjectures the following stronger result.

\begin{conjecture}
Let $\mathbf{P}$ be an arbitrary totally-sewn neighbourly $4$-polytope in ${\mathbb{E}}^4$, and $\mathbf{o}$ be an arbitrary interior point of $\mathbf{P}$. Then there exist $9$ hyperplanes of ${\mathbb{E}}^4$ such that each face of $\mathbf{P}$, can be strictly separated from $\mathbf{o}$ by at least one of the $9$ hyperplanes.
\end{conjecture}

\subsection{General results} \label{sec:general}
Before we state results on the illumination number of convex bodies in ${\mathbb{E}}^{d}$, we take a little detour. We need Rogers' estimate  \cite{rogers} of the infimum $\theta(\mathbf{K})$ of the covering density of ${\mathbb{E}}^{d}$ by translates of the convex body $\mathbf{K}$, namely, for $d\geq 2$, \footnote{The bound on $\theta(\mathbf{K})$ has been improved to $\theta(\mathbf{K})\leq d\ln d + d\ln\ln d + d + o(d)$ by G. Fejes T\'{o}th \cite{fejestoth}.}
\[\theta(\mathbf{K})\leq d(\ln d + \ln\ln d + 5) 
\]
and the Rogers--Shephard inequality \cite{rogers-shephard} 
\[\vol_{d}(\mathbf{K}-\mathbf{K})\leq \binom{2d}{d}\vol_{d}(\mathbf{K})
\]
on the $d$-dimensional volume $\vol_{d}(\cdot)$ of the difference body $\mathbf{K}-\mathbf{K}$ of $\mathbf{K}$.  

It was rather a coincidence, at least from the point of view of the Illumination Conjecture, when in 1964 Erd\H{o}s and Rogers \cite{erdos2} proved the following theorem. In order to state their theorem in a proper form we need to introduce the following notion. If we are given a covering of a space by a system of sets, the {\it star number} of the covering is the supremum, over sets of the system, of the cardinals of the numbers of sets of the system meeting a set of the system (see \cite{erdos2}). On the one hand, the standard Lebesgue €œbrick-laying€ construction provides an example, for each positive integer $d$, of a lattice covering of ${\mathbb{E}}^{d}$ by closed cubes with star number $2^{d+1} - 1$. On the other hand, Theorem 1 of \cite{erdos2} states that the star number of a lattice covering of ${\mathbb{E}}^{d}$ by translates of a centrally symmetric convex body is always at least $2^{d+1} - 1$. However, from our point of view, the main result of \cite{erdos2} is the one under Theorem 2 which (combined with some observations from \cite{erdos1} and the Rogers--Shephard inequality \cite{rogers}) reads as follows. 

\begin{theorem}\label{rogers} Let $\mathbf{K}$ be a convex body in the $d$-dimensional Euclidean space ${\mathbb{E}}^{d}$, $d\geq 2$. Then there exists a covering of ${\mathbb{E}}^{d}$ by translates of $\mathbf{K}$ with star number at most 
\[
\frac{\vol_{d} (\mathbf{K} - \mathbf{K})}{\vol_{d}(\mathbf{K})} (d \ln d+d\ln\ln d+5d+1) \leq {2d\choose d} (d\ln d+d\ln\ln d+5d+1).
\]
Moreover, for sufficiently large $d$, $5d$ can be replaced by $4d$.
\end{theorem}

The periodic and probabilistic construction on which Theorem \ref{rogers} is based gives also the following. 

\begin{corollary}\label{rogerscor} If $\mathbf{K}$ is an arbitrary convex body in ${\mathbb{E}}^{d}$, $d \geq 2$, then
\begin{equation}\label{ill-bound}
I(\mathbf{K}) \leq \frac{\vol_{d} (\mathbf{K} - \mathbf{K})}{\vol_{d}(\mathbf{K})} d(\ln d+\ln\ln d+5)\leq {2d\choose d}d(\ln d+\ln\ln d+5)= O(4^{d}\sqrt{d}\ln d).
\end{equation}
Moreover, for sufficiently large $d$, $5d$ can be replaced by $4d$.
\end{corollary}

Note that the bound given in Corollary \ref{rogerscor} can also be obtained from the more general result of Rogers and Zong \cite{rogers-zong}, which states that for $d$-dimensional convex bodies $\mathbf{K}$ and $\mathbf{L}$, $d\geq 2$, one can cover $\mathbf{K}$ by $N(\mathbf{K},\mathbf{L})$ translates \footnote{$N(\mathbf{K},\mathbf{L})$ is called the {\it covering number of $\mathbf{K}$ by $\mathbf{L}$}.} of $\mathbf{L}$ such that 
\[N(\mathbf{K},\mathbf{L})\leq \frac{\vol_{d}(\mathbf{K} - \mathbf{L})}{\vol_{d}(\mathbf{L})}{\theta}(\mathbf{L}). 
\]

For the sake of completeness we also mention the inequality 
$$I(\mathbf{K}) \leq (d + 1)d^{d-1} - (d - 1)(d - 2)^{d -1}$$
due to Lassak \cite{lassak3}, which is valid for an arbitrary convex body $\mathbf{K}$ in ${\mathbb{E}}^{d}$, $d \geq 2$, and€™ is (somewhat) better than the estimate of Corollary \ref{rogerscor} for some small values of $d$.

Since, for a centrally symmetric convex body $\mathbf{K}$, $\frac{{\textnormal{vol}}(\mathbf{K}-\mathbf{K})}{\vol_{d}(\mathbf{K})} = 2^d$, we have the following improved upper bound on the illumination number of such convex bodies.  

\begin{corollary}\label{rogerscor2}
If $\mathbf{K}$ is a centrally symmetric convex body in ${\mathbb{E}}^{d}$, $d \geq 2$, then
\begin{equation}\label{ill-symmetric}
I(\mathbf{K}) \leq \frac{\vol_{d} (\mathbf{K} - \mathbf{K})}{\vol_{d}(\mathbf{K})}d(\ln d+\ln\ln d+5)= 2^{d}d (\ln d+\ln\ln d+5)= O(2^{d}d\ln d).
\end{equation}
\end{corollary}

The above upper bound is fairly close to the conjectured value of $2^{d}$. However, most convex bodies are far from being symmetric and so, in general, one may wonder whether the Illumination Conjecture is true at all, especially for large $d$. Thus, it was important progress when Schramm \cite{schramm} managed to prove the Illumination Conjecture for all convex bodies of constant width in all dimensions at least 16. In fact, he proved the following inequality. 

\begin{theorem}\label{constantwidth} If $\mathbf{W}$ is an arbitrary convex body of constant width in ${\mathbb{E}}^{d}$ , $d \geq 3$, then 
\[
I(\mathbf{W})\leq 5d\sqrt{d}(4+\ln d) \left(\frac{3}{2}\right)^{\frac{d}{2}}.
\]
\end{theorem}

By taking a closer look of Schramm's elegant paper \cite{schramm} and making the necessary modifications, the first author \cite{Be12f} somewhat improved the upper bound of Theorem~\ref{constantwidth}, but more importantly he succeeded in extending that estimate to the following family of convex bodies (called the family of {\it fat spindle convex bodies}) that is much larger than the family of convex bodies of constant width. Thus, we have the following generalization of Theorem~\ref{constantwidth} proved in \cite{Be12f}.

\begin{theorem}\label{19-2}
Let $X\subset\mathbb{E}^{d}$, $d\ge 3$ be an arbitrary compact set with ${\rm diam}(X)\le 1$ and let $\mathbf{B}[X]$ be the intersection of the closed $d$-dimensional unit balls centered at the points of $X$. Then $$I(\mathbf{B}[X])< 4\left(\frac{\pi}{3}\right)^{\frac{1}{2}}d^{\frac{3}{2}}(3+\ln d)\left(\frac{3}{2}\right)^{\frac{d}{2}}<5d^{\frac{3}{2}}(4+\ln d)\left(\frac{3}{2}\right)^{\frac{d}{2}}.$$
\end{theorem}

On the one hand, $4\left(\frac{\pi}{3}\right)^{\frac{1}{2}}d^{\frac{3}{2}}(3+\ln d)\left(\frac{3}{2}\right)^{\frac{d}{2}}<2^d$ for all $d\ge 15$. (Moreover, for every $\epsilon>0$ if $d$ is sufficiently large, then $I(\mathbf{B}[X])<\left(\sqrt{1.5}+\epsilon\right)^d=(1.224\ldots +\epsilon)^d$.) On the other hand, based on the elegant construction of Kahn and Kalai \cite{kahn}, it is known (see \cite{AiZi}), that if $d$ is sufficiently large, then there exists a finite subset $X''$ of $\{0,1\}^d$ in $\mathbb{E}^{d}$ such that any partition of $X''$ into parts of smaller diameter requires more than $(1.2)^{\sqrt{d}}$ parts. Let $X'$ be the (positive) homothetic copy of $X''$ having unit diameter and let $X$ be the (not necessarily unique) convex body of constant width one containing $X'$. Then it follows via standard arguments that $I(\mathbf{B}[X])>(1.2)^{\sqrt{d}}$ with $X=\mathbf{B}[X]$.

Recall that a convex polytope is called a \textit{belt polytope} if to each side of any of its 2-faces there exists a parallel (opposite) side on the same 2-face. This class of polytopes is wider than the class of zonotopes. Moreover, it is easy to see that any convex body of ${\mathbb{E}}^{d}$ can be represented as a limit of a covergent sequence of belt polytopes with respect to the Hausdorff metric in ${\mathbb{E}}^{d}$. The following theorem on belt polytopes was proved by Martini in \cite{martini2}. The result that it extends to the class of convex bodies, called belt bodies (including zonoids), is due to Boltyanski \cite{boltyanski-add, boltyanski5, boltyanski2}. (See also \cite{boltyanski7} for a somewhat sharper result on the illumination numbers of belt bodies.)

\begin{theorem} \label{belt} Let $\mathbf{P}$ be an arbitrary $d$-dimensional belt polytope (resp., belt body) different from a parallelotope in ${\mathbb{E}}^{d}$, $d \geq 2$. Then
\[
I(\mathbf{P}) \leq 3 \cdot 2^{d-2}.
\]
\end{theorem}

Now, let $\mathbf{K}$ be an arbitrary convex body in ${\mathbb{E}}^{d}$ and let ${\cal T}(\mathbf{K})$ be the family of all translates of $\mathbf{K}$ in ${\mathbb{E}}^{d}$. The {\it Helly dimension} ${\textnormal{him}}(\mathbf{K})$ of $\mathbf{K}$ (\cite{soltan0}) is the smallest integer $h$ such that for any finite family ${\cal{F}} \subseteq {\cal T}(\mathbf{K})$ with cardinality greater than $h + 1$ the following assertion holds: if every $h + 1$ members of ${\cal{F}}$ have a point in common, then all the members of ${\cal{F}}$ have a point in common. As is well known $1 \leq {\textnormal{him}}(\mathbf{K}) \leq d$. Using this notion Boltyanski \cite{boltyanski6} gave a proof of the following theorem.

\begin{theorem}\label{him} Let $\mathbf{K}$ be a convex body with ${\rm him}(\mathbf{K}) = 2$ in ${\mathbb{E}}^{d}$, $d \geq 3$. Then 
\[
I(\mathbf{K}) \leq 2^{d} - 2^{d-2}.
\]
\end{theorem}

In fact, in \cite{boltyanski6} Boltyanski conjectures the following more general inequality. 

\begin{conjecture} \label{himconjecture} Let $\mathbf{K}$ be a convex body with ${\textnormal{him}}(\mathbf{K}) = h > 2$ in ${\mathbb{E}}^{d}$, $d \geq 3$. Then
\[
I(\mathbf{K}) \leq 2^{d} - 2^{d-h}.
\]
\end{conjecture}

The first author and Bisztriczky gave a proof of the Illumination Conjecture for the class of dual cyclic polytopes in \cite{bezdek-bisz}. Their upper bound for the illumination numbers of dual cyclic polytopes has been improved by Talata in \cite{talata}. So, we have the following statement.

\begin{theorem} \label{dual} The illumination number of any $d$-dimensional dual cyclic polytope is at most $\frac{(d+1)^2}{2}$, for all $d \geq 2$.
\end{theorem}

%In connection with the results of this section quite a number of questions remain open including the following ones.

%\begin{problem}
%\item{(a)} What are the illumination numbers of cyclic polytopes?
%\item{(b)} Can one give a proof of the separation conjecture for zonotopes (resp., belt polytopes)?
%\item{(c)} Is there a way to prove the separation conjecture for $0/1$-polytopes?
%\end{problem}

\section{On some relatives of the illumination number}\label{sec:relatives}

\subsection{llumination by affine subspaces}

Let $\mathbf{K}$ be a convex body in $\mathbb{E}^{d}$, $d\geq 2$. The following definitions were introduced by the first named author in \cite{B94} (see also \cite{bezdek-conj1} that introduced the concept of the first definition below).

 Let $L\subset \mathbb{E}^{d}\setminus \mathbf{K}$ be an affine subspace of dimension $l$, $0\le l\le d-1$. Then $L$ illuminates the boundary point $\mathbf{q}$ of $\mathbf{K}$ if there exists a point $\mathbf{p}$ of $L$ that illuminates $\mathbf{q}$ on the boundary of  $\mathbf{K}$. Moreover, we say that the affine subspaces $L_1, L_2, \dots , L_n$ of dimension $l$ with $L_i\subset \mathbb{E}^{d}\setminus \mathbf{K}, 1\le i\le n$ illuminate $\mathbf{K}$ if every boundary point of $\mathbf{K}$ is illuminated by at least one of the affine subspaces $L_1, L_2, \dots , L_n$. Finally, let $I_l(\mathbf{K})$ be the smallest positive integer $n$ for which there exist $n$ affine subspaces of dimension $l$ say, $L_1, L_2, \dots , L_n$ such that $L_i\subset \mathbb{E}^{d}\setminus \mathbf{K}$ for all $1\le i\le n$ and  $L_1, L_2, \dots , L_n$ illuminate $\mathbf{K}$. Then $I_l(\mathbf{K})$ is called the {\sl $l$-dimensional illumination number} of $\mathbf{K}$ and the sequence $I_0(\mathbf{K}), I_1(\mathbf{K}), \dots , I_{d-2}(\mathbf{K}), I_{d-1}(\mathbf{K})$ is called the {\it successive illumination numbers} of $\mathbf{K}$. Obviously, $I(\mathbf{K})=I_0(\mathbf{K})\ge I_1(\mathbf{K})\ge \dots \ge I_{d-2}(\mathbf{K})\ge I_{d-1}(\mathbf{K})=2$.

  Recall that $\mathbb{S}^{d-1}$ denotes the unit sphere centered at the origin of $\mathbb{E}^{d}$. Let ${HS}^l\subset \mathbb{S}^{d-1}$ be an $l$-dimensional open great-hemisphere of $\mathbb{S}^{d-1}$, where $0\le l\le d-1$. Then ${HS}^l$ illuminates the boundary point $\mathbf{q}$ of $\mathbf{K}$ if there exists a unit vector $\mathbf{v}\in {HS}^l$ that illuminates $\mathbf{q}$, in other words, for which it is true that the halfline emanating from $\mathbf{q}$ and having direction vector $\mathbf{v}$ intersects the interior of $\mathbf{K}$. Moreover, we say that the $l$-dimensional open great-hemispheres $ {HS}^l_1, {HS}^l_2, \dots , {HS}^l_n$ of $\mathbb{S}^{d-1}$ illuminate $\mathbf{K}$ if each boundary point of $\mathbf{K}$ is illuminated by at least one of the open great-hemispheres $ {HS}^l_1, {HS}^l_2, \dots , {HS}^l_n$. Finally, let $I'_l(\mathbf{K})$ be the smallest number of $l$-dimensional open great-hemispheres of $\mathbb{S}^{d-1}$ that illuminate $\mathbf{K}$. Obviously, $I'_0(\mathbf{K})\ge I'_1(\mathbf{K})\ge \dots \ge I'_{d-2}(\mathbf{K})\ge I'_{d-1}(\mathbf{K})=2$.

Let $L\subset \mathbb{E}^{d}$ be a linear subspace of dimension $l$, $0\le l\le d-1$ in  $\mathbb{E}^{d}$. The {\it $l$-codimensional circumscribed  cylinder} of $\mathbf{K}$ generated by $L$ is the union of translates of $L$ that have a nonempty intersection with $\mathbf{K}$. Then let $C_l(\mathbf{K})$ be the smallest number of translates of the interiors of some $l$-codimensional circumscribed  cylinders of $\mathbf{K}$ the union of which contains $\mathbf{K}$. Obviously, $C_0(\mathbf{K})\ge C_1(\mathbf{K})\ge \dots \ge C_{d-2}(\mathbf{K})\ge C_{d-1}(\mathbf{K})=2$.

The following theorem, which was proved in \cite{B94}, collects the basic information known about the quantities just introduced. 

\begin{theorem} \label{21}
Let $\mathbf{K}$ be an arbitrary  convex body of $\mathbb{E}^{d}$. Then
\item(i) $I_l(\mathbf{K})= I'_l(\mathbf{K})= C_l(\mathbf{K})$, for all $0\le l\le d-1$.
\item(ii) $\lceil{\frac{d+1}{l+1}}\rceil\le I_l(\mathbf{K})$, for all $0\le l\le d-1$, with equality for any smooth $\mathbf{K}$.
\item(iii) $I_{d-2}(\mathbf{K})=2$, for all $d\ge 3$.
\end{theorem}

The Generalized Illumination Conjecture was phrased by the first named author in \cite{B94} as follows.

\begin{conjecture}[{\bf Generalized Illumination Conjecture}]\label{Bezdek-general-illumination-conjecture}  
Let $\mathbf{K}$ be an arbitrary  convex body and $\mathbf{C}^d$ be a $d$-di\-men\-si\-on\-al affine cube in $\mathbb{E}^{d}$. Then
$$I_l(\mathbf{K})\le I_l(\mathbf{C}^d)$$
holds for all $0\le l\le d-1$.
\end{conjecture}

 The above conjecture was proved for zonotopes and zonoids in \cite{B94}. The results of parts (i) and (ii) of the next theorem are taken from \cite{B94}, where they were proved for zonotopes (resp., zonoids).  However, in the light of the more recent works in \cite{boltyanski5} and \cite{boltyanski7} these results extend to the class of belt polytopes (resp., belt bodies) in a rather straightforward way so we present them in that form. The lower bound of part (iii) was proved in \cite{B94} and the upper bound of part (iii) is the major result of \cite{K}. Finally, part (iv) was proved in \cite{BK93}.

\begin{theorem}
Let $\mathbf{M}$ be a belt polytope (resp., belt body) and $\mathbf{C}^d$ be a $d$-dimensional affine cube in $\mathbb{E}^{d}$. Then

\item(i) $I_l(\mathbf{M})\le I_l(\mathbf{C}^d)$  holds for all $0\le l\le d-1$.

\item(ii) $I_{\lfloor\frac{d}{2}\rfloor}(\mathbf{M})=\dots =I_{d-1}(\mathbf{M})=2$.

\item(iii) $\frac{2^d}{\sum_{i=0}^{l}{\binom{d}{i}}}\le I_l(\mathbf{C}^d)\le K(d, l)$, where $K(d,l)$ denotes the minimum cardinality of binary codes of length $d$ with covering radius $l$, $0\le l\le d-1$.

\item(iv) $I_1(\mathbf{C}^d)=\frac{2^d}{d+1}$, provided that $d+1=2^m$.
\end{theorem}

\subsection{`X-raying' the problem}\label{sec:x-ray relative}

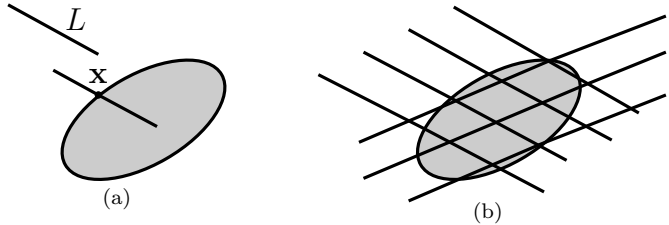
\begin{figure}[htb]
	\label{fig:lines}
	\centering
\begin{tabular}{ccc}
\subfigure [ ]{
\begin{tikzpicture}[scale = 0.6, baseline]
\large
\draw[fill=black!20, rotate=30, very thick] (0,0) ellipse [x radius=2cm, y radius=1cm];
\draw[fill=black] (-1,0.55) circle [radius=2pt];
\draw[very thick] (-2,1.1) -- (-1,0.55) node[above] {$\mathbf{x}$} -- (0.3,-.15); 
%\draw[->, very thick] (-2,2.05) -- (-0.7,1.5) node[midway,above] {$\vec{\ell}$}; 
%\draw[->, very thick] (-2,2) -- (-0.2,1);  
\draw[very thick] (-3, 2.55) -- (-2,2) -- (-1,1.45) node[midway,above] {$L$};  
\normalsize
\end{tikzpicture}
} &    &
\subfigure [ ]{
\begin{tikzpicture}[scale = 0.6, baseline]
\large
\draw[fill=black!20, rotate=30, very thick] (0,0) ellipse [x radius=2cm, y radius=1cm];
\draw[very thick] (-1,2.5) -- (1.1,1.3) -- (3.1,0.15); 
\draw[very thick] (-2,2) -- (-0.2,1) -- (2.5,-.45); 
\draw[very thick] (-3,1.5) -- (-1.1,0.5) -- (1.5,-.9); 
\draw[very thick] (-4,1) -- (-1.65,-0.2) -- (1,-1.6); 

\draw[very thick] (3.7,2.3) -- (1.1,1.3) -- (-3.1,-0.5); 
\draw[very thick] (3.7,1.5) -- (1.85,.75) -- (-3,-1.3); 
\draw[very thick] (3.7,.55) -- (0.9,-.55) -- (-2,-1.8); 

%\draw[->, very thick] (-3.5,-2) -- (-1.65,-0.2); 
%\draw[->, very thick] (-2.3,-2) -- (-1.6,-1.3); 
%\draw[->, very thick] (-1.1,-2) -- (-.5,-1.3); 
%\draw[->, very thick] (-.2,-2) -- (0.9,-.7); 
\normalsize
\end{tikzpicture}
}
\end{tabular}
\caption{(a) X-raying a boundary point $\mathbf{x}$ of $\mathbf{K}$ along a line $L$, (b) $X(\mathbf{K}) = 2$.}
\end{figure}

In 1972, the \textit{X-ray number} of convex bodies was introduced by P. Soltan as follows (see \cite{martini1}). Let $\mathbf{K}$ be a convex body of ${\mathbb{E}}^{d}$, $d \geq 2$, and $L \subset {\mathbb{E}}^{d}$ be a line through the origin of ${\mathbb{E}}^{d}$. We say that the boundary point $\mathbf{x} \in \mathbf{K}$ is X-rayed along $L$ if the line parallel to $L$ passing through $\mathbf{x}$ intersects the interior of $\mathbf{K}$. The X-ray number $X(\mathbf{K})$ of $\mathbf{K}$ is the smallest number of lines such that every boundary point of $\mathbf{K}$ is X-rayed along at least one of these lines. Clearly, $X(\mathbf{K}) \geq d$. Moreover, it is easy to see that this bound is attained by any smooth convex body. On the other hand, if ${\mathbf{C}}^{d}$ is a $d$-dimensional (affine) cube and $F$ is one of its $(d - 2)$-dimensional faces, then the X-ray number of ${\mathbf{C}}^{d}\setminus F$, the convex hull of the set of vertices of ${\mathbf{C}}^{d}$ that do not belong to $F$, is $3\cdot 2^{d-ˆ'2}$.

In 1994, the first author and Zamfirescu \cite{bezdek-z} published the following conjecture.

\begin{conjecture}[{\bf X-ray Conjecture}]\label{x-ray} The X-ray number of any convex body in ${\mathbb{E}}^{d}$ is at most $3 \cdot 2^{d-2}$. 
\end{conjecture}

The X-ray Conjecture is proved only in the plane and it is open in higher dimensions. Here we note that the inequalities 
$$X(\mathbf{K}) \leq I(\mathbf{K}) \leq 2X(\mathbf{K})$$ 
hold for any convex body $\mathbf{K} \subset {\mathbb{E}}^{d}$. In other words, any proper progress on the X-ray Conjecture would imply progress on the Illumination Conjecture and vice versa. We also note that a natural way to prove the X-ray Conjecture would be to show that any convex body $\mathbf{K} \subset {\mathbb{E}}^{d}$ can be illuminated by $3 \cdot 2^{d-2}$ pairs of pairwise opposite directions. 

The main results of \cite{bezdek-kiss} on the X-ray number can be summarized as follows. In order to state them properly we need to recall two basic notions. Let $\mathbf{K}$ be a convex body in ${\mathbb{E}}^{d}$  and let $F$ be a face of $\mathbf{K}$.The {\it Gauss image} $\nu (F)$ of the face $F$ is the set of all points (i.e., unit vectors) $\mathbf{u}$ of the $(d - 1)$-dimensional unit sphere ${\mathbb{S}}^{d-1} \subset {\mathbb{E}}^{d}$ centered at the origin $\mathbf{o}$ of ${\mathbb{E}}^d$ for which the supporting hyperplane of $K$ with outer normal vector $\mathbf{u}$ contains $F$. It is easy to see that the Gauss images of distinct faces of $\mathbf{K}$ have disjoint relative interiors in ${\mathbb{S}}^{d-1}$ and $\nu(F)$ is compact and spherically convex for any face $F$. Let $C \subset {\mathbb{S}}^{d-1}$ be a set of finitely many points. Then the {\it covering radius} of $C$ is the smallest positive real number $r$ with the property that the family of spherical balls of radii $r$ centered at the points of $C$ covers ${\mathbb{S}}^{d-1}$. 

\begin{theorem} \label{gaussimage} Let $\mathbf{K}\subset {\mathbb{E}}^{d}$, $d\geq 3$, be a convex body and let $r$ be a positive real number with the property that the Gauss image $\nu(F)$ of any face $F$ of $\mathbf{K}$ can be covered by a spherical ball of radius $r$ in ${\mathbb{S}}^{d-1}$. Moreover, assume that there exist $2m$ pairwise antipodal points of $\mathbb{S}^{d-1}$ with covering radius $R$ satisfying the inequality $r+R\le\frac{\pi}{2}$. Then $X(\mathbf{K})\le m$. In particular,
if there are $2m$
pairwise antipodal points on $\mathbb{S}^{d-1}$ with covering radius
$R$ satisfying the inequality $R\leq \pi /2-r_{d-1}$, where $r_{d-1}=\arccos \sqrt{\frac{d+1}{2d}}$ is the circumradius of a
regular $(d-1)$-dimensional spherical simplex of edge length $\pi /3$, then $X(\mathbf{W})\le m$ holds for any
convex body $\mathbf{W}$ of constant width in $\mathbb{E}^d$.
\end{theorem}

\begin{theorem} \label{x-ray2} If $\mathbf{W}$ is an arbitrary convex body of constant width in ${\mathbb{E}}^{3}$, then $X(\mathbf{W}) = 3$. If $\mathbf{W}$ is any convex body of constant width in ${\mathbb{E}}^{4}$, then $4 \leq X(\mathbf{W}) \leq 6$. Moreover, if $\mathbf{W}$ is a convex body of constant width in ${\mathbb{E}}^{d}$ with $d=5,6$, then $d\leq X(\mathbf{W})\leq 2^{d-ˆ'1}$.
\end{theorem}

\begin{corollary} If $\mathbf{W}$ is an arbitrary convex body of constant width in ${\mathbb{E}}^{3}$, then $4 \leq I(\mathbf{W}) \leq6$. If $\mathbf{W}$ is any convex body of constant width in ${\mathbb{E}}^{4}$, then $5 \leq I(\mathbf{W}) \leq 12$. Moreover, if $\mathbf{W}$ is a convex body of constant width in ${\mathbb{E}}^{d}$ with $d=5,6$, then $d+1\leq I(\mathbf{W})\leq 2^{d}$.
\end{corollary}

It would be interesting to extend the method described in the paper \cite{bezdek-kiss} for the next few dimensions (more exactly, for the dimensions $7 \leq d \leq 14$) in particular, because in these dimensions neither the X-ray Conjecture nor the Illumination Conjecture is known to hold for convex bodies of constant width. 

From the proof of Theorem~\ref{Lassak-illumination} it follows in a straightforward way that if $\mathbf{K}$ is a centrally symmetric convex body in ${\mathbb{E}}^3$, then $X(\mathbf{K})\leq 4$. On the other hand, very recently Trelford \cite{Tr14} proved the following related result.
\begin{theorem}
If $\mathbf{K}$ is a convex body symmetric about a plane in ${\mathbb{E}}^3$, then $X(\mathbf{K})\leq 6$.
\end{theorem}

\subsection{Other relatives}\label{sec:other relatives}
\subsubsection{$t$-covering and $t$-illumination numbers}\label{sec:t}
In Section \ref{sec:intro}, we found that the least number of smaller positive homothets of a convex body $\mathbf{K}$ required to cover it equals the minimum number of translates of the interior of $\mathbf{K}$ needed to cover $\mathbf{K}$. Is this number also equal to the the minimum number $t(\mathbf{K})$ of translates of $\mathbf{K}$ that are different from $\mathbf{K}$ and are needed to cover $\mathbf{K}$? 

Despite being a very natural question, the problem of economical translative coverings have not attracted much attention. To our knowledge, the first systematic study of these was carried out quite recently by Lassak, Martini and Spirova \cite{lassak-translative} who called them $t$-coverings and also introduced the corresponding illumination concept, called $t$-illumination, as follows: A boundary point $\mathbf{x}$ of a convex body $\mathbf{K}$ of ${\mathbb{E}}^{d}$ is {\it $t$-illuminated} by a direction $\mathbf{v}\in\mathbb{S}^{d-1}$ if there exists a different point $\mathbf{y} \in \mathbf{K}$ such that the vector $\mathbf{y}-\mathbf{x}$ has the same direction as $\mathbf{v}$ (i.e., $\mathbf{y}-\mathbf{x}=\lambda\mathbf{v}$, for some $\lambda>0$). The minimum number $i(\mathbf{K})$ of directions needed to $t$-illuminate the entire boundary of $\mathbf{K}$ is called its $t$-illumination number. The connection between $t$-covering and $t$-illumination is summarized in the next result \cite{lassak-translative}. Note that a convex body $\mathbf{K}$ is said to be {\it strictly convex} if for any two points of $\mathbf{K}$ the open line segment connecting them belongs to the interior of $\mathbf{K}$. 

\begin{theorem}
\item{(i)} If $\mathbf{K}$ is a planar convex body, then $i(\mathbf{K}) = t(\mathbf{K})$. 

\item{(ii)} If $\mathbf{K}$ is a $d$-dimensional strictly convex body, $d\geq 3$, then $i(\mathbf{K})= t(\mathbf{K})$. 

\item{(iii)} If $\mathbf{K}$ is a $d$-dimensional convex body, $d\geq 3$, then $i(\mathbf{K})\leq t(\mathbf{K})$, where the equality does not hold in general.   
\end{theorem}

Clearly, $t(\mathbf{K})\leq I(\mathbf{K})$. In the same paper \cite{lassak-translative}, the following results were obtained about the relationship between $I(\cdot)$ and $t(\cdot)$. 

\begin{theorem}
\item{(i)} If $\mathbf{K}$ is a planar convex body, then $t(\mathbf{K}) = I(\mathbf{K})$ if and only if $\mathbf{K}$ contains no parallel boundary segments. 

\item{(i)} If $\mathbf{K}$ is a strictly convex body of ${\mathbb{E}}^{d}$, $d\geq3$, then $t(\mathbf{K})= I(\mathbf{K})$.

\item{(iii)} If $\mathbf{K}$ is a convex body of ${\mathbb{E}}^{d}$, $d\geq3$ that does not have parallel boundary segments, then $t(\mathbf{K})= I(\mathbf{K})$. 
 \end{theorem}

However, in general the following remains unanswered \cite{lassak-translative}. 

\begin{problem}
Characterize the convex bodies $\mathbf{K}$ for which $t(\mathbf{K})=I(\mathbf{K})$. 
\end{problem}

In \cite{martini-t-illumination}, the notion of $t$-illumination was refined into $t$-central illumination and strict $t$-illumination and the corresponding illumination numbers were defined. The paper also introduced metric versions of the classical, $t$-central and strict $t$-illumination numbers and investigated their properties at length. The interested reader is referred to \cite{martini-t-illumination} for details.  

\subsubsection{Blocking numbers}\label{sec:block}
The {\it blocking number} $\beta(\mathbf{K})$ \cite{zong-block} of a convex body $\mathbf{K}$ is defined as the minimum number of nonoverlapping translates of $\mathbf{K}$ that can be brought into contact with the boundary of $K$ so as to block any other translate of $\mathbf{K}$ from touching $\mathbf{K}$. Since $\beta(\mathbf{K})=\beta(\mathbf{K}-\mathbf{K})$ and $\mathbf{K}-\mathbf{K}$ is $\mathbf{o}$-symmetric, it suffices to consider the blocking numbers of $\mathbf{o}$-symmetric convex bodies only.   

For any $\mathbf{o}$-symmetric convex body $\mathbf{K}$, the relation $I(\mathbf{K})\leq \beta(\mathbf{K})$ holds \cite{zong-block}, while no such relationship exists for general convex bodies. Zong \cite{zong-block} conjectured the following. 

\begin{conjecture}
For any $d$-dimensional convex body $\mathbf{K}$, 
$$2d\leq \beta(\mathbf{K})\leq 2^{d},$$
and $\beta(\mathbf{K})=2^{d}$ if and only if $\mathbf{K}$ is a $d$-dimensional cube.  
\end{conjecture}

If true, Zong's Conjecture would imply the Illumination Conjecture for $\mathbf{o}$-symmetric convex bodies. Some of the known values of the blocking number include $\beta(\mathbf{K})=2^{d}$, if $\mathbf{K}$ is a $d$-dimensional cube; $\beta(\mathbf{K})=6$, if $\mathbf{K}$ is a 3-dimensional ball; and $\beta(\mathbf{K})=9$, if $\mathbf{K}$ is a 4-dimensional ball \cite{dalla-zong}. Some other values and estimates are obtained in \cite{yu-blocking}. 

Several generalizations of the blocking number have been proposed. The smallest number of non-overlapping translates of $\mathbf{K}$ such that the interior of $\mathbf{K}$ is disjoint from the interiors of the translates and they can block any other translate from touching $\mathbf{K}$ is denoted by $\beta_{1}(\mathbf{K})$; the smallest number of translates all of which touch $\mathbf{K}$ at its boundary such that they can block any other translate from touching $\mathbf{K}$ is denoted by $\beta_{2}(\mathbf{K})$; whereas, $\beta_{3}(\mathbf{K})$ denotes the smallest number of translates all of which are non-overlapping with $\mathbf{K}$ such that they can block any other translate from touching $\mathbf{K}$ \cite{yu-zong}. If in the original definition of blocking number, translates are replaced by homothets with homothety ratio $\alpha >0$ we get the {\it generalized blocking number} $\beta^{\alpha}(\mathbf{K})$ \cite{boroczky-zong}, and if we allow the homothets to overlap, we get the {\it generalized $\alpha$-blocking number} $\beta^{\alpha}_{2}(\mathbf{K})$ \cite{wu-blocking}. 

Recently, Wu \cite{wu-blocking} showed that if $\mathbf{K}$ and $\mathbf{L}$ are $\mathbf{o}$-symmetric convex bodies that are sufficiently close to each other in the Banach--Mazur sense\footnote{\label{bmh} See relation \eqref{banach-mazur} and the discussion preceding it in Section \ref{sec:parameters} for an introduction to the Banach--Mazur distance of convex bodies. Note that Wu uses the Hausdorff distance between convex bodies to state his result. However, it can be shown that $\mathbf{K}$ and $\mathbf{L}$ are close to each other in the Banach--Mazur sense if and only if there exist affine images of them that are close in the Hausdorff sense. Since the illumination and blocking numbers are affine invariants, we can restate Wu's results in the language of Banach--Mazur distance.} then there exists $\alpha >0$ (depending on $\mathbf{K}$) such that 
\[I(\mathbf{K})\leq \beta^{\alpha}_{2}(\mathbf{L}). 
\]

This gives a series of upper bounds on the illumination number of symmetric convex bodies and a possible way to circumvent the lack of lower semicontinuity of $I(\cdot)$ (see Section \ref{sec:parameters} for a discussion of the continuity of the illumination number).

\subsubsection{Fractional covering and illumination}\label{sec:fractional}
Nasz\'{o}di \cite{naszodi-fractional} introduced the fractional illumination number and Arstein-Avidan with Raz \cite{artstein-avidan1} and with Slomka \cite{artstein-avidan2} introduced weighted covering numbers. Both formalisms can be used to study a fractional analogue of the illumination problem. In fact, the Fractional Illumination Conjecture for $\mathbf{o}$-symmetric convex bodies was proved in \cite{naszodi-fractional}, while the case of equality was characterized in \cite{artstein-avidan2}. We omit the details as it would lead to a lengthy diversion from the main subject matter.

\section{Quantifying illumination and covering}\label{sec:quantify}
\subsection{The illumination and covering parameters}\label{sec:parameters}
It can be seen that in the definition of illumination number $I(\mathbf{K})$, the distance of light sources from $\mathbf{K}$ plays no role whatsoever. Starting with a relatively small number of light sources, it makes sense to quantify how far they need to be from $\mathbf{K}$ in order to illuminate it. This is the idea behind the {\it illumination parameter} as defined by the first author \cite{bezdek-illumination1}.

Let $\mathbf{K}$ be an $\mathbf{o}$-symmetric convex body. Then the norm of $x\in {\mathbb{E}}^{d}$ generated by $K$ is defined as 
\[\left\|\mathbf{x}\right\|_{\mathbf{K}}=\inf \{\lambda>0: \mathbf{x}\in \lambda \mathbf{K}\}
\]
and provides a good estimate of how far a point $\mathbf{x}$ is from $\mathbf{K}$. 

The illumination parameter $\ill(\mathbf{K})$ of an $\mathbf{o}$-symmetric convex body $\mathbf{K}$ estimates how well $\mathbf{K}$ can be illuminated by relatively few point sources lying as close to $\mathbf{K}$ on average as possible. 
\[
\ill(\mathbf{K})=\inf \left\{\sum_{i}\left\|\mathbf{p_{i}}\right\|_{K} : \{\mathbf{p_{i}}\}\textnormal{ illuminates } \mathbf{K}, \mathbf{p_{i}}\in {\mathbb{E}}^{d}\right\}, 
\]  

Clearly, $I(\mathbf{K})\leq\ill(\mathbf{K})$ holds for any $\mathbf{o}$-symmetric convex body $\mathbf{K}$. In the papers \cite{bezdek-illumination2,kiss1}, the illumination parameters of $\mathbf{o}$-symmetric Platonic solids have been determined. In \cite{bezdek-illumination1} a tight upper bound was obtained for the illumination parameter of planar $\mathbf{o}$-symmetric convex bodies.  

\begin{theorem}
If $\mathbf{K}$ is an $\mathbf{o}$-symmetric planar convex body, then $\ill(\mathbf{K}) \leq 6$ with equality for any affine regular convex hexagon. 
\end{theorem} 

The corresponding problem in dimension 3 and higher is wide open. The following conjecture is due to Kiss and de Wet \cite{kiss1}.   

\begin{conjecture}
The illumination parameter of any $\mathbf{o}$-symmetric 3-dimensional convex body is at most 12. 
\end{conjecture} 

However, for smooth $\mathbf{o}$-symmetric convex bodies in any dimension $d\geq 2$, the first named author and Litvak \cite{bezdek-litvak} found an upper bound, which was later improved to the following asymptotically sharp bound by Gluskin and Litvak \cite{gluskin-litvak}. 

\begin{theorem}\label{vein}
For any smooth $\mathbf{o}$-symmetric $d$-dimensional convex body $\mathbf{K}$, 
\[\ill(\mathbf{K})\leq 24 d^{3/2}.
\]
\end{theorem}
  
Translating the above quantification ideas from illumination into the setting of covering, Swanepoel \cite{swanepoel1} introduced the {\it covering parameter} of a convex body as follows. 
\[C(\mathbf{K})=\inf \left\{\sum_{i}(1-\lambda_{i})^{-1}:\mathbf{K}\subseteq \bigcup_{i}(\lambda_{i}\mathbf{K}+\mathbf{t_{i}}), 0<\lambda_{i}<1,\mathbf{t_{i}}\in {\mathbb{E}}^{d}\right\}. 
\]  

Thus large homothets are penalized in the same way as the far off light sources are penalized in the definition of illumination parameter. Note that here $\mathbf{K}$ need not be $\mathbf{o}$-symmetric. In the same paper, Swanepoel obtained the following Rogers-type upper bounds on $C(\mathbf{K})$ when $d\geq 2$. 

\begin{theorem} 
\begin{equation}\label{swanepoel1} 
C(\mathbf{K})<\left\{\begin{split} e2^{d}d(d+1)(\ln d+\ln \ln d + 5)=O(2^{d}d^{2}\ln d), \ \ \ \ \ \ \ \ \ & \textnormal{ if } \mathbf{K} \textnormal{ is } o\textnormal{-symmetric},\\
e\binom{2d}{d}d(d+1)(\ln d+\ln \ln d + 5)=O(4^{d}d^{3/2}\ln d), \ \ & \ \textnormal{otherwise}.
\end{split} \right.\ \ \ 
\end{equation}
\end{theorem}

He further showed that if $K$ is $o$-symmetric, then
\begin{equation}\label{swanepoel2}
\ill(\mathbf{K}) \leq 2 C(\mathbf{K}), 
\end{equation} 
and therefore, $\ill(\mathbf{K})=O(2^{d}d^{2}\ln d)$.  

Based on the above results, it is natural to study the following quantitative analogue of the illumination conjecture that was proposed by Swanepoel \cite{swanepoel1}. 

\begin{conjecture}[{\bf Quantitative Illumination Conjecture}]
For any $\mathbf{o}$-symmetric $d$-dimensional convex body $\mathbf{K}$, $\ill(\mathbf{K}) = O(2^{d})$.
\end{conjecture}

%\subsection{Properties of $I(\cdot)$, $\ill(\cdot)$ and $C(\cdot)$}\label{sec:properties}
Before proceeding further, we introduce some terminology and notations. Let us use ${\cal{K}}^{d}$ and ${\cal{C}}^{d}$ respectively to denote the set of all $d$-dimensional convex bodies and the set of all such bodies that are $\mathbf{o}$-symmetric. In this section, we consider some of the important properties of the illumination number and the covering parameter as functionals defined on ${\cal{K}}^{d}$ and the illumination parameter as a functional on ${\cal{C}}^{d}$. The first observation is that the three quantities are affine invariants (as are several other quantities dealing with the covering and illumination of convex bodies). That is, if $A:{\mathbb{E}}^{d}\to {\mathbb{E}}^{d}$ is an affine transformation and $\mathbf{K}$ is any $d$-dimensional convex body, then $I(\mathbf{K})=I(A(\mathbf{K}))$, $\ill(\mathbf{K})=\ill(A(\mathbf{K}))$ and $C(\mathbf{K})=C(A(\mathbf{K}))$. 

Due to this affine invariance, whenever we refer to a convex body $\mathbf{K}$, whatever we say about the covering and illumination of $\mathbf{K}$ is true for all affine images of $\mathbf{K}$. In the sequel, $\mathbf{B}^{d}$ denotes a $d$-dimensional unit ball\footnote{Without loss of generality, we can assume $\mathbf{B}^{d}$ to be a unit ball centred at the origin. In what follows, we use the symbol $\mathbf{B}^{d}$ to denote a $d$-dimensional unit ball as well as its affine images called ellipsoids.}, $\mathbf{C}^{d}$ a $d$-dimensional cube and $\ell$ a line segment (which is a convex body in ${\cal{K}}^{1}$) up to an affine transformation. 

The Banach--Mazur distance $d_{BM}$ provides a multiplicative metric\footnote{One can turn the Banach--Mazur distance into an additive metric by applying $\ln(\cdot)$. However, we make no attempt to do that.} on ${\cal{K}}^{d}$  and is used to study the continuity properties of affine invariant functionals on ${\cal{K}}^{d}$. For $\mathbf{K}, \mathbf{L}\in {\cal{K}}^{d}$, it is given by 
\begin{equation}\label{banach-mazur} 
d_{BM}(\mathbf{K}, \mathbf{L})= \inf \left\{\delta \geq 1 : \mathbf{L}-\mathbf{b}\subseteq T(\mathbf{K}-\mathbf{a})\subseteq \delta (\mathbf{L}-\mathbf{b}), \mathbf{a} \in \mathbf{K}, \mathbf{b}\in \mathbf{L}\right\},
\end{equation}
where the infimum is taken over all invertible linear operators $T:{\mathbb{E}}^{d}\to {\mathbb{E}}^{d}$ \cite[Page 589]{schneider1}. 

In the remainder of this paper, ${\cal{K}}^{d}$ (resp., ${\cal{C}}^{d}$) is considered as a metric space under the Banach--Mazur distance. Since continuity of a functional can provide valuable insight into its behaviour, it is of considerable interest to check the continuity of $I(\cdot)$, $\ill(\cdot)$ and $C(\cdot)$. Unfortunately, by Example \ref{discontinuity}, the first two quantities are known to be discontinuous, while nothing is known about the continuity of the third. 

\begin{example}[Smoothed cubes and spiky balls]\label{discontinuity}
In ${\cal{K}}^{d}$, consider a sequence $(\mathbf{C_{n}})_{n\in {\mathbb{N}}}$ of `smoothed' $d$-dimensional cubes that approaches $\mathbf{C}^d$ in the Banach--Mazur sense. Since the smoothed cubes are smooth convex bodies, all the terms of the sequence have illumination number $d+1$. However, $I(\mathbf{C}^{d})=2^{d}$, which shows that $I(\cdot)$ is not continuous. 

Recently, Nasz\'{o}di \cite{naszodi-spiky} constructed a class of $d$-dimensional $\mathbf{o}$-symmetric bodies, that he refers to as `spiky balls'. Pick $N$ points $ \mathbf{x_{1}},\ldots, \mathbf{x_{N}}$ independently and uniformly with respect to the Haar probability measure on the $(d-1)$-dimensional unit sphere $\mathbb{S}^{d-1}$ centred at the origin $\mathbf{o}$. Then a spiky ball corresponding to a real number $D>1$ is defined as 
\[\mathbf{K}=\conv\left(\{\pm \mathbf{x_{i}} : i=1,\ldots, N\} \cup \frac{1}{D}\mathbf{B}^{d}\right).
\]

Straightaway we observe that $K$ is $\mathbf{o}$-symmetric and satisfies  $d_{BM}(\mathbf{K},\mathbf{B}^{d})<D$. Nasz\'{o}di showed that $I(\mathbf{K})\geq c^{d}$, where $c>1$ is a constant depending on $d$ and $D$. Thus we have a sequence of spiky balls approaching $\mathbf{B}^{d}$ in Banach--Mazur distance such that each spiky ball has an exponential illumination number. Since by Theorem \ref{vein}, $\ill(\mathbf{B}^{d})=O(d^{3/2})$ and $\ill(\mathbf{K})\geq I(\mathbf{K})$ we see that $\ill(\cdot)$ is not continuous.  
\end{example}

We can state the following about the continuity of $I(\cdot)$ \cite{boltyanski2}.  

\begin{theorem}
The functional $I(\cdot)$ is upper semicontinuous\footnote{Again, the original statement of this result is in terms of Hausdorff distance. However, based on the discussion in footnote \ref{bmh}, we can use the Banach--Mazur distance instead.} on ${\mathbb{E}}^{d}$, for all $d\geq 2$. 
\end{theorem} 

Despite the usefulness of the covering parameter, not much is known about it. For instance, we do not know whether $C(\cdot)$ is lower or upper semicontinuous on ${\cal{K}}^{d}$ and the only known exact value is $C(\mathbf{C}^{d})=2^{d+1}$. Thus there is a need to propose a more refined quantitative version of homothetic covering for convex bodies. Section \ref{sec:index} describes how we address this need.

\subsection{The covering index}\label{sec:index}
The concepts and results presented in this section appear in our recent paper \cite{bezdek-khan1}. As stated at the end of Section \ref{sec:parameters}, the aim here is to come up with a more refined quantification of covering in terms of the covering index with the Covering Conjecture as the eventual goal. The {\it covering index} of a convex body $\mathbf{K}$ in $\Eu^d$ combines the notions of the covering parameter $C(\mathbf{K})$ and the $m$-covering number $\gamma_{m}(\mathbf{K})$ under the unusual, but highly useful, constraint $\gamma_{m}(\mathbf{K})\leq 1/2$, where $\gamma_{m}(\mathbf{K})=\inf \left\{\lambda >0: \mathbf{K}\subseteq \bigcup_{i=1}^{m}(\lambda \mathbf{K}+\mathbf{t_{i}}), \mathbf{t_{i}}\in {\mathbb{E}}^{d}, i=1,\ldots, m\right\}$ is the smallest positive homothety ratio needed to cover $\mathbf{K}$ by $m$ positive homothets. (See Section \ref{sec:computational} for a detailed discussion of $\gamma_{m}(\cdot)$.)

\begin{definition}[Covering index]\label{coin-def}
Let $\mathbf{K}$ be a $d$-dimensional convex body. We write $N_{\lambda}(\mathbf{K})$ to denote the covering number $N(\mathbf{K},\lambda \mathbf{K})$, for any $0< \lambda \leq 1$. We define the \textit{covering index} of $\mathbf{K}$ as 
\begin{align*}
\coin(\mathbf{K})&=\inf \left\{\frac{m}{1-\gamma_{m}(\mathbf{K})}: \gamma_{m}(\mathbf{K})\leq 1/2, m\in \mathbb{N}\right\}\\
&= \inf \left\{\frac{N_{\lambda}(\mathbf{K})}{1-\lambda}: 0<\lambda\leq 1/2\right\}. 
\end{align*}
\end{definition}

Intuitively, $\coin(\mathbf{K})$ measures how $\mathbf{K}$ can be covered by a relatively small number of positive homothets all corresponding to the same relatively small homothety ratio. The reader may be a bit surprised to see the restriction $\gamma_{m}(\mathbf{K})\leq 1/2$. In \cite{bezdek-khan1}, it was observed that if we start with $\gamma_{m}(\mathbf{K})\leq \lambda <1$, for some $\lambda$ close to 1 in the definition of $\coin(\mathbf{K})$ and then decrease $\lambda$, the properties of the resulting quantity significantly change when $\lambda =1/2$, at which point we can say a lot about the continuity and maximum and minimum values of the quantity. It was also observed that decreasing $\lambda$ further does not change these characteristics. Thus $1/2$ can be thought of as a threshold at which the characteristics of the covering problem change. 

Note that for $\mathbf{K}\in {\cal{C}}^{d}$,  
\[ I(\mathbf{K})\leq \ill(\mathbf{K})\leq 2C(\mathbf{K}) \leq 2\coin(\mathbf{K}), \footnote{In fact, one can obtain $\ill(\mathbf{K})\leq \frac{3}{2}\coin(\mathbf{K})$ by suitably modifying the proof of \eqref{swanepoel2} which appears in \cite{swanepoel1}.}\]
\noindent and in general for $\mathbf{K}\in {\cal{K}}^{d}$, 
\[ I(\mathbf{K})\leq C(\mathbf{K}) \leq \coin(\mathbf{K}). \]

The following result shows that a lot can be said about the Banach--Mazur continuity of $\coin(\cdot)$. Based on this, $\coin(\cdot)$ seems to be the `nicest' of all the functionals of covering and illumination of convex bodies discussed here. 

\begin{theorem}\label{continuity} Let $d$ be any positive integer.  
\item(i) Define $I_{\mathbf{K}}=\{i: \gamma_{i}(\mathbf{K})\leq 1/2\}=\{i: \mathbf{K}\in {\cal{K}}^{d}_{i}\}$, for any $d$-dimensional convex body $\mathbf{K}$. If $I_{\mathbf{L}}\subseteq I_{\mathbf{K}}$, for some $\mathbf{K}, \mathbf{L} \in {\cal{K}}^{d}$, then $\coin(\mathbf{K}) \leq \frac{2d_{BM}(\mathbf{K},\mathbf{L})-1}{d_{BM}(\mathbf{K},\mathbf{L})} \coin(\mathbf{L}) \leq d_{BM}(\mathbf{K},\mathbf{L}) \coin(\mathbf{L})$. 
\item(ii) The functional $\coin:{\cal{K}}^{d}\to \mathbb{R}$ is lower semicontinuous for all $d$.  
\item(iii) Define ${\cal{K}}^{d*}:=\left\{\mathbf{K}\in {\cal{K}}^{d}: \gamma_{m}(\mathbf{K})\neq 1/2, m\in {\mathbb{N}} \right\}$. Then the functional $\coin:{\cal{K}}^{d*}\to \mathbb{R}$ is continuous for all $d$.
\end{theorem}

We now present some results showing that $\coin(\cdot)$ behaves very nicely with forming direct sums, Minkowski sums and cylinders of convex bodies, making it possible to compute the exact values and estimates of $\coin(\cdot)$ for higher dimensional convex bodies from the covering indices of lower dimensional convex bodies. 

\begin{theorem}\label{productnew}   
\item(i) Let ${\mathbb{E}}^{d}={\mathbb{L}}_{1} \oplus \cdots \oplus {\mathbb{L}}_{n}$ be a decomposition of ${\mathbb{E}}^{d}$ into the direct vector sum of its linear subspaces ${\mathbb{L}}_{i}$ and let $\mathbf{K_{i}}\subseteq {\mathbb{L}}_{i}$ be convex bodies such that $\Gamma=\max\{ \gamma_{m_{i}}(\mathbf{K_{i}}):1\leq i\leq n\}$. Then 
\begin{equation}\label{eq:product1}
\begin{split}
\max_{1\leq i\leq n} \{\coin(\mathbf{K_{i}})\} \leq  \coin(\mathbf{K_{1}}\oplus \cdots \oplus \mathbf{K_{n}}) \leq \inf_{\lambda \leq \frac{1}{2}} \frac{\prod_{i=1}^{n} N_{\lambda}(\mathbf{K_{i}})}{1-\lambda}  \leq \frac{\prod_{i=1}^{n} N_{\Gamma}(\mathbf{K_{i}})}{1-\Gamma} < \prod_{i=1}^{n} \coin(\mathbf{K_{i}}). 
\end{split}
\end{equation}

\item(ii) The first two upper bounds in (\ref{eq:product1}) are tight. Moreover, the second inequality in \eqref{eq:product1} becomes an equality if any $n-1$ of the $\mathbf{K_{i}}$'s are tightly covered

\item(iii) Recall that $\ell\in {\cal{K}}^{1}$ denotes a line segment. If $\mathbf{K}$ is any convex body, then $\coin(\mathbf{K}\oplus \ell) = 4N_{1/2}(\mathbf{K})$.

\item(iv) Let the convex body $K$ be the Minkowski sum of the convex bodies $\mathbf{K_{1}}, \ldots , \mathbf{K_{n}} \in {\cal{K}}^{d}$ and $\Gamma$ be as in part (i). Then 
\begin{equation}\label{mink}
\coin(\mathbf{K})  \leq \inf_{\lambda \leq \frac{1}{2}} \frac{\prod_{i=1}^{n} N_{\lambda}(\mathbf{K_{i}})}{1-\lambda} \leq \frac{\prod_{i=1}^{n} N_{\Gamma}(\mathbf{K_{i}})}{1-\Gamma} < \prod_{i=1}^{n} \coin(\mathbf{K_{i}}). 
\end{equation} 
\end{theorem}

The notion of {\it tightly covered convex bodies} introduced in \cite{bezdek-khan1} plays a critical role in Theorem \ref{productnew} (ii)-(iii).  

\begin{definition}\label{tight}
We say that a convex body $\mathbf{K}\in {\cal{K}}^{d}$ is \textit{tightly covered} if for any $0<\lambda <1$, $\mathbf{K}$ contains $N_{\lambda}(\mathbf{K})$ points no two of which belong to the same homothet of $\mathbf{K}$ with homothety ratio $\lambda$. 
\end{definition}

In \cite{bezdek-khan1}, it was noted that not all convex bodies are tightly covered (e.g., $\mathbf{B}^{2}$ is not), $\ell\in {\cal{K}}^{1}$ is tightly covered and so is the $d$-dimensional cube $\mathbf{C}^{d}$, for any $d\geq 2$. Do other examples exist? 

\begin{problem}\label{other-tight}
For some $d\geq 2$, find a tightly covered convex body $\mathbf{K}\in {\cal{K}}^{d}$ other than $\mathbf{C}^{d}$ or show that no such convex body exists.  
\end{problem}

Since $\coin$ is a lower semicontinuous functional defined on the compact space ${\cal{K}}^{d}$, it is guaranteed to achieve its infimum over ${\cal{K}}^{d}$. It turns out that in addition to determining minimizers in all dimensions, we can also find a maximizer in the planar case. 
  
\begin{theorem}\label{cube}
\item(i) Let $d$ be any positive integer and $\mathbf{K}\in {\cal{K}}^{d}$. Then $\coin(\mathbf{C}^{d})=2^{d+1}\leq \coin(\mathbf{K})$ and thus $d$-cubes minimize the covering index in all dimensions. 

\item(ii) If $\mathbf{K}$ is a planar convex body then $\coin(\mathbf{K})\leq \coin(\mathbf{B}^{2})=14$.
\end{theorem}

Since $\mathbf{B}^{2}$ maximizes the covering index in the plane, it can be asked if the same is true for $\mathbf{B}^{d}$ in higher dimensions. 

\begin{problem}\label{ball-max}
For any $d$-dimensional convex body $\mathbf{K}$, prove or disprove that $\coin(\mathbf{K})\leq \coin(\mathbf{B}^{d})$ holds. 
\end{problem} 

Since $\coin(\mathbf{B}^{d})=O(2^{d}d^{3/2}\ln d)$ \cite{bezdek-khan1}, a positive answer to Problem \ref{ball-max} would considerably improve the best known upper bound on the illumination number $I(\mathbf{K})=O(4^{d}\sqrt{d} \ln d)$ when $\mathbf{K}$ is a general $d$-dimensional convex body to within a factor $\sqrt{d}$ of the bound $I(\mathbf{K})=O(2^{d}d \ln d)$ when $\mathbf{K}$ is $\mathbf{o}$-symmetric. This gives us a way to closing in on the Illumination Conjecture for general convex bodies. 

If we replace the restriction $\gamma_{m}(\mathbf{K})\leq 1/2$ from the definition of the covering index with the more usual $\gamma_{m}(\mathbf{K})<1$, the resulting quantity is called the {\it weak covering index}, denoted by $\wcoin(\mathbf{K})$. As the name suggests, the weak covering index loses some of the most important properties of the covering index. For instance, no suitable analogue of Theorem \ref{productnew} (iii) exists for $\wcoin(\cdot)$. As a result, we can only estimate the weak covering index of cylinders. Also the discussed aspects of continuity of the covering index seem to be lost for the weak covering index. Last, but not the least, unlike the covering index we cannot say much at all about the maximizers and minimizers of the weak covering index. 

In the end, we would like to mention that fractional analogues of the covering index and the weak covering index were introduced in \cite{bezdek-khan2}. Just like fractional illumination number, we do not discuss these here due to limitation of space.

\subsection{Cylindrical covering parameters}\label{sec:kth-CC}
So far in Section \ref{sec:parameters}-\ref{sec:index}, we have discussed some quantitative versions of the illumination number. The aim of this section is to introduce a quantification of the X-ray number. This quantification has the added advantage of connecting the X-ray problem with the Tarski's plank problem and its relatives (see \cite[Chapter 4]{bezdek-book}).   
 
Given a linear subspace $E\subseteq {\Eu}^d$ we denote the orthogonal projection on $E$ by $P_E$ and the orthogonal complement of $E$ by $E^{\perp}$.
Given $0<k<d$, define a {\it $k$-codimensional cylinder} $\mathbf{C}$ as a set, which can be
presented in the form $\mathbf{C} = B + H$, where $H$ is a $k$-dimensional linear subspace
of ${\Eu}^{d}$ and  $B$ is a measurable set in $E: =H^{\perp}$.  Given a convex
body $\KK$ and a $k$-co\-di\-men\-si\-o\-nal cylinder $\mathbf{C}= B + H$ denote
the cross-sectional volume of $\mathbf{C}$ with respect to $\KK$ by
$$
 \crv _{\KK} (\mathbf{C}) := \frac{\vol _{d-k} (\mathbf{C}\cap E)}{\vol _{d-k} (P_E
 \KK)} =\frac{\vol _{d-k} (P_E \mathbf{C})}{\vol _{d-k} (P_E
 \KK)} = \frac{\vol _{d-k} (B)}{\vol _{d-k} (P_E \KK)} .
$$

We note that if $T:\Eu^d\to\Eu^d$ is an invertible affine map, then $ \crv _{\KK} (\mathbf{C})= \crv _{T(\KK)} (T(\mathbf{C}))$. Now we introduce the following.  

\begin{definition}[$k$-th Cylindrical Covering Parameter]\label{k-Xindex}
Let $0<k<d$ and $\KK$ be a convex body in $\Eu^d$. Then the $k$-th cylindrical covering parameter of $\KK$ is labelled by $\cyl_k(\KK)$ and it is defined as follows:
$$\cyl_k(\KK)=\inf_{\bigcup_{i=1}^{n}\mathbf{C}_i}\bigg\{\sum_{i=1}^n \crv _{\KK} (\mathbf{C}_i)\ :\ \KK\subseteq\bigcup_{i=1}^{n}\mathbf{C}_i, \ \mathbf{C}_i\ {\rm is\ a\ }k-{\rm codimensional}\  {\rm cylinder},\ i=1,\ldots,n\bigg\}.$$
\end{definition}

We observe that if $T:\Eu^d\to\Eu^d$ is an invertible affine map, then $\cyl_k(\KK)=\cyl_k(T(\KK))$. Furthermore, it is clear that $\cyl_k(\KK)\leq 1$ holds for any convex body $\KK$ in $\Eu^d$ and for any $0<k<d$. In terms of X-raying, one can think of $\cyl_{k}(\KK)$ as the minimum of the `sum of sizes' of $(d-k)$-dimensional X-raying windows needed to X-ray $\KK$. 

%To be more precise, consider a collection $\{\mathbf{C}_{i}=B_{i}+H_{i}:i=1,\ldots, n\}$ of $k$-codimensional cylinders in $\Eu^{d}$, where $H_{i}$ is a $k$-dimensional linear subspace of $\Eu^d$ and $B_{i}$ is a measurable set in $H_{i}^\perp$, $i=1,\ldots, n$. Each $\mathbf{C}_{i}$ is a union of lines (1-dimensional affine subspaces) in $\Eu^d$ parallel to $H_{i}$ and orthogonal to $H_{i}^\perp$. These lines can X-ray $\KK$ and thus we can think of $\mathbf{C}_{i}$ as an X-ray beam through a window $B_{i}$ of size $\crv_\KK(\mathbf{C}_{i})$. In addition, it is clear that $\KK$ is X-rayed by these beams if and only if the interiors $\inter(\mathbf{C}_{i})$ of the $\mathbf{C}_i$'s cover $\KK$. Finally, note that if $\KK\subseteq\bigcup_{i=1}^{n}\mathbf{C}_i$, then there exist open cylinders $\mathbf{D}_{i}$, $i=1,\ldots, n$, such that  $\sum_{i=1}^n \crv _{\KK} (\mathbf{C}_i) = \sum_{i=1}^n \crv _{\KK} (\mathbf{D}_i)$. 

Recall that a $(d-1)$-codimensional cylinder of $\Eu^d$ is also called a {\it plank} for the reason that it is the set of points lying between two parallel hyperplanes in $\Eu^d$. The width of a plank is simply the distance between the two parallel hyperplanes. In a remarkable paper \cite{Ba51}, Bang has given an elegant proof of the Plank Conjecture of Tarski showing that if a convex body is covered by finitely many planks in $\Eu^d$, then the sum of the widths of the planks is at least as large as the minimal width of the body, which is the smallest distance between two parallel supporting hyperplanes of the given convex body. A celebrated extension of Bang's theorem to $d$-dimensional normed spaces has been given by Ball in \cite{Bal91}. In his paper \cite{Ba51}, Bang raises his so-called Affine Plank Conjecture, which in terms of our notation can be phrased as follows. 

\begin{conjecture}[{\bf Affine Plank Conjecture}]\label{Bang-conjecture}
If $\KK$ is a convex body in $\Eu^d$, then $\cyl_{d-1}(\KK)=1$.
\end{conjecture}

Now, Ball's celebrated plank theorem (\cite{Bal91}) can be stated as follows.

\begin{theorem}
If $\KK$ is an $\mathbf{o}$-symmetric convex body in $\Eu^d$, then $\cyl_{d-1}(\KK)=1$.
\end{theorem}

Bang \cite{Ba51} also raised the important related question of whether the sum of the base areas of finitely many (1-codimensional) cylinders covering a 3-dimensional convex body is at least half of the minimum area of a 2-dimensional projection of the body. This, in terms of our terminology, reads as follows.

\begin{conjecture}[{\bf $1$-Codimensional Cylinder Covering Conjecture}]\label{Bang-cylinder-covering-conjecture}
If $\KK$ is a convex body in $\Eu^3$, then $\cyl_{1}(\KK)\geq\frac{1}{2}$.
\end{conjecture}

If true, then Bang's estimate is sharp due to a covering of a regular tetrahedron by two cylinders described in \cite{Ba51}. In connection with Conjecture~\ref{Bang-cylinder-covering-conjecture} the first named author and Litvak have proved the following general estimates in \cite{BeLi09}.

\begin{theorem}\label{Bezdek-Litvak-cylinder-covering}
Let $0<k<d$ and $\KK$ be a convex body in $\Eu^d$. Then $\cyl_k(\KK)\geq \frac{1}{{d \choose k}}$. 
\end{theorem}

Furthermore, it is proved in \cite{BeLi09} that if $\KK$ is an ellipsoid in $\Eu^d$, then $\cyl_1(\KK)=1$. 
Akopyan, Karasev and Petrov (\cite{AKP}) have recently proved that if $\KK$ is an ellipsoid in $\Eu^d$, then $\cyl_2(\KK)=1$. They have put forward: 

\begin{conjecture}[{\bf Ellipsoid Conjecture}]\label{ellipsoid-conjecture}
If $\KK$ is an ellipsoid in $\Eu^d$, then $\cyl_k(\KK)=1$ for all $2<k<d$.
\end{conjecture}

\section{A computer-based approach} \label{sec:computational} 
Given a positive integer $m$, Lassak \cite{lassak-gamma} introduced the \textit{$m$-covering number} of a convex body $\mathbf{K}$ as the minimal positive homothety ratio needed to cover $\mathbf{K}$ by $m$ positive homothets. That is, 
\[\gamma_{m}(\mathbf{K})=\inf \left\{\lambda >0: \mathbf{K}\subseteq \bigcup_{i=1}^{m}(\lambda \mathbf{K}+\mathbf{t_{i}}), \mathbf{t_{i}}\in {\mathbb{E}}^{d}, i=1,\ldots, m\right\}.
\]
Lassak showed that the $m$-covering number is well-defined and studied the special case $m=4$ for planar convex bodies. It should be noted that special values of this quantity had been considered by several authors in the past. For instance, in the 70's and 80's the first named author showed that $\gamma_{5}(\mathbf{B}^{2})=0.609382\ldots$\footnote{{\it Cover the Spot} is a popular carnival game in the United States. The objective is to cover a given circular spot by 5 circular disks of smaller radius. It seems that by determining $\gamma_{5}(\mathbf{B}^{2})$, the first named author was unwittingly providing the optimal solution for {\it Cover the Spot}!} and $\gamma_{6}(\mathbf{B}^{2})=0.555905\ldots$ \cite{bezdek1,bezdek2}. 

\vspace{-6mm}
\begin{figure}[htb]
	\centering
		\includegraphics[scale=0.4]{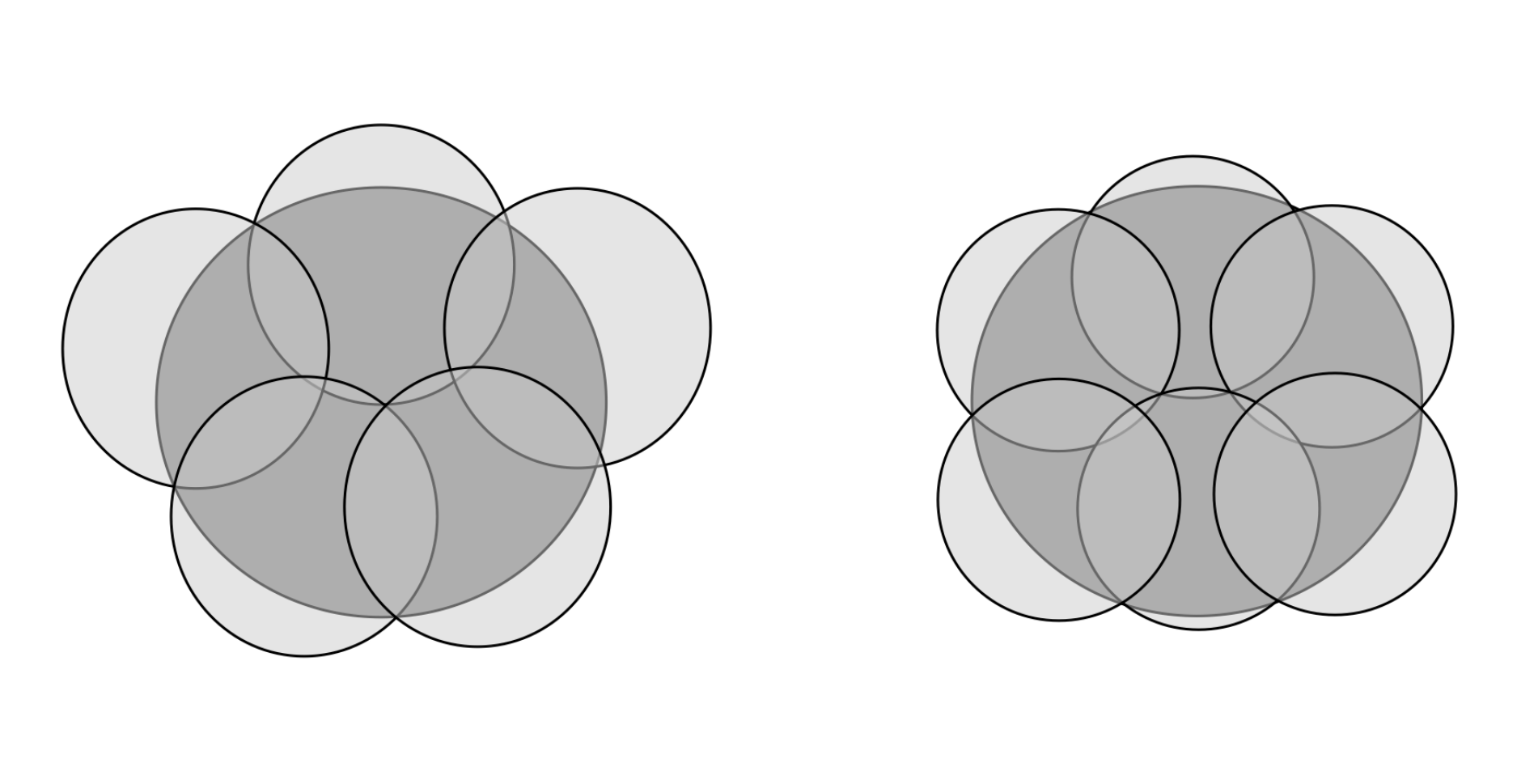}
	\label{fig:circles}
\vspace{-6mm}
	\caption{Optimal configurations that demonstrate $\gamma_{5}(\mathbf{B}^{2}) = 0.609382\ldots$ and $\gamma_{6}(\mathbf{B}^{2}) = 0.555905\ldots$}
\end{figure}

Zong \cite{zong1} studied $\gamma_{m}:{\cal{K}}^{d} \to \mathbb{R}$ as a functional and proved it to be uniformly continuous for all $m$ and $d$. He did not use the term $m$-covering number for $\gamma_{m}(\mathbf{K})$ and simply referred to it as the smallest positive homothety ratio. In \cite{bezdek-khan1}, we proved the following stronger result. 

\begin{theorem} \label{aux}
For any $K,L\in {\cal{K}}^{d}$, $\gamma_{m}(K)\leq d_{BM}(K,L)\gamma_{m}(L)$ holds and so $\gamma_{m}$ is Lipschitz continuous on ${\cal{K}}^{d}$ with $\frac{d^2-1}{2\ln d}$ as a Lipschitz constant and 
$$\left|\gamma_{m}(K)-\gamma_{m}(L)\right|\le d_{BM}(K,L)-1\le \frac{d^2-1}{2\ln d}\ln \left(d_{BM}(K,L)\right),$$
for all $K, L \in {\cal{K}}^{d}$. 
\end{theorem}

Further properties and some variants of $\gamma_{m}(\cdot)$ are discussed in the recent papers \cite{he, wu-chinese}. For instance, it has been shown in \cite{wu-chinese} that for any $d$-dimensional convex polytope $\mathbf{P}$ with $m$ vertices, we have 
\[\gamma_{m}(\mathbf{K})\leq \frac{d-1}{d}. 
\] 

Obviously, any $\mathbf{K} \in {\cal{K}}^{d}$ can be covered by $2^{d}$ smaller positive homothets if and only if $\gamma_{2^{d}}(\mathbf{K})<1$. Zong used these ideas to propose a possible computer-based approach to attack the Covering Conjecture \cite{zong1}. 

Recall that in a metric space, such as ${\cal{K}}^{d}$, an {\it $\epsilon$-net} $\xi$ is a finite or infinite subset of ${\cal{K}}^{d}$ such that the union of closed balls of radius $\epsilon$ centered at elements of $\xi$ covers the whole space. Thus if an $\epsilon$-net exists, any element of ${\cal{K}}^{d}$ is within Banach--Mazur distance $\epsilon$ of some element of the cover. The key idea of the procedure proposed by Zong is the construction of a finite $\epsilon$-net of ${\cal{K}}^{d}$ whose elements are convex polytopes, for every real number\footnote{Recall that the Banach--Mazur distance is a multiplicative metric and so the condition $\epsilon>0$ is replaced by the equivalent $\epsilon > 1$ condition.} $\epsilon > 1$ and positive integer $d$. Here we describe the construction briefly. 

\begin{figure}[htb]
	\centering
		\includegraphics[scale=0.3]{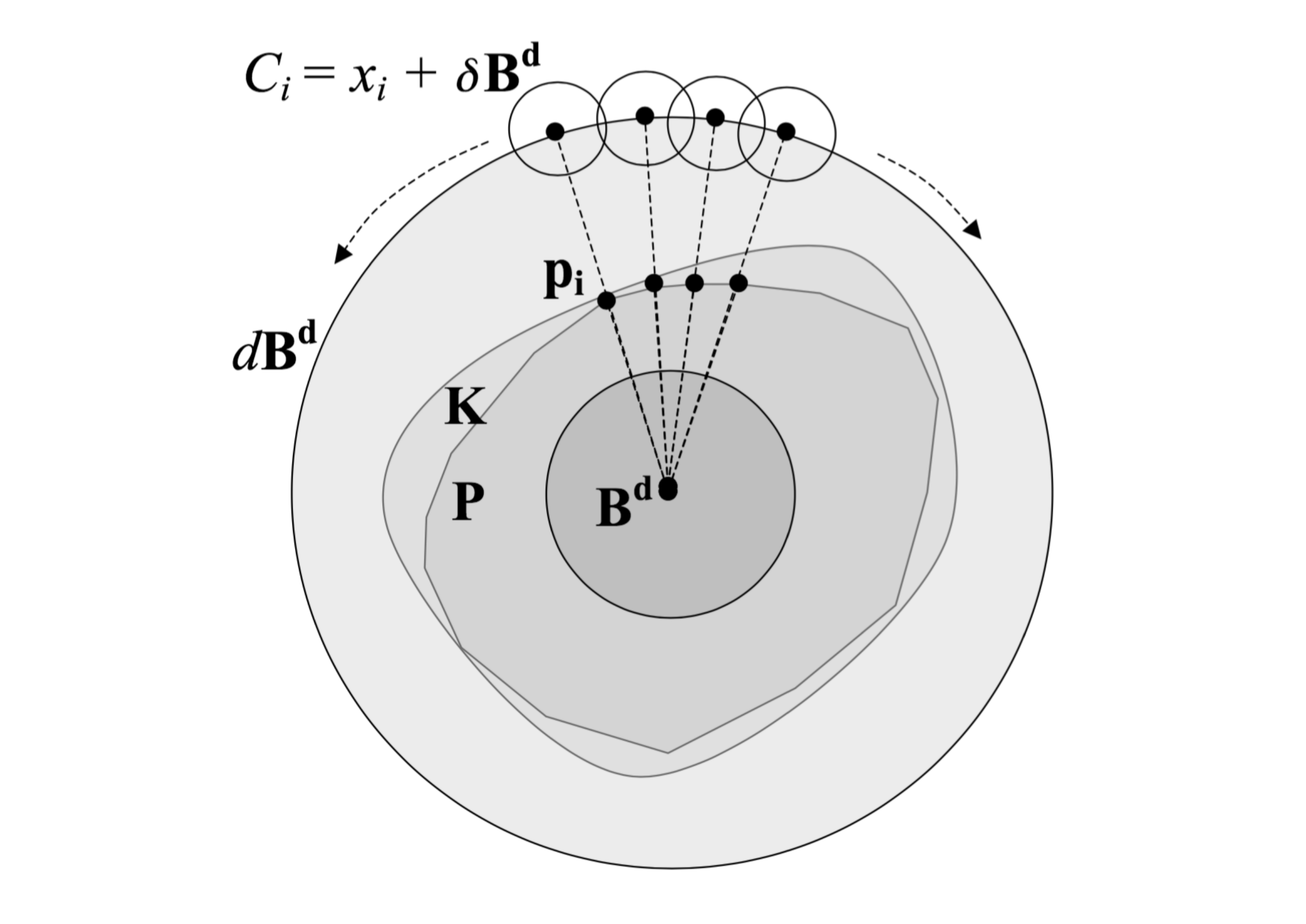}
	\label{fig:zong}
\vspace{-2mm}	
	\caption{Construction of an $\epsilon$-net of ${\cal{K}}^{d}$}
\end{figure}

We first take an affine image of a $d$-dimensional convex body $\mathbf{K}$ that is sandwiched between the unit ball $\mathbf{B^{d}}$ centered at the origin and the ball $d\mathbf{B^{d}}$ with radius $d$. Such an image always exists by John's ellipsoid theorem. Then we take a covering $\{C_{1},\ldots,C_{m}\}$ of the boundary of $d\mathbf{B^{d}}$ with spherical caps $C_{i}$ as shown in Figure \ref{fig:zong}. The centers of the caps $C_{i}$ are joined to the origin by lines $\{L_{1},\ldots, L_{m}\}$ and a large number of equidistant points are taken on the lines $L_{i}$. We denote by $\mathbf{p_{i}}$ the point lying in $\mathbf{K}\cap L_i$ that is farthest from the origin. Then the convex hull $\mathbf{P}=\conv\{\mathbf{p_{i}}\}$ is the required element of our $\epsilon$-net. Zong \cite{zong1} showed that by taking $m$ large enough and increasing the number of points on $L_i$ we can ensure $d_{BM}(\mathbf{K},\mathbf{P})\leq \epsilon$.  

He then notes that if we manage to construct a finite $\epsilon$-net $\xi=\{\mathbf{P_{i}}:i=1,\ldots, j\}$ of ${\cal{K}}^{d}$, satisfying $\gamma_{2^{d}}(\mathbf{P_{i}})\leq c_d$ for some $c_{d}<1$ and sufficiently small $\epsilon$, then $\gamma_{2^{d}}(\mathbf{K})<1$ would hold for every $\mathbf{K}\in {\cal{K}}^{d}$. This would imply that the Covering Conjecture is true in dimension $d$. 

The following is a four step approach suggested by Zong \cite{zong1}.  

\vspace{4mm}
\noindent \textbf{Zong's Program:} 
\begin{enumerate}
\item For a given dimension such as $d = 3$, investigate (with the assistance of a computer) $\gamma_{2^{d}}(\mathbf{K})$ for some particular convex bodies $\mathbf{K}$ and choose a candidate constant $c_{d}$. 
\item Choose a suitable $\epsilon$. 
\item Construct an $\epsilon$-net $\xi$ of sufficiently small cardinality. 
\item Check (with the assistance of a computer) that the minimal $\gamma_{2^{d}}$-value over all elements of $\xi$ is bounded above by $c_{d}$. 
\end{enumerate}

Indeed this approach appears to be promising and, to the authors' knowledge, is a first attempt at a computer-based resolution of the Covering Conjecture. However, Zong's program is not without its pitfalls. For one, it would take an extensive computational experiment to come up with a good candidate constant $c_{d}$. Secondly, Zong's $\epsilon$-net construction leads to a net with exponentially large number of elements. In fact, using B\"or\"oczky and Wintsche's estimate \cite{boroczky1} on the number of caps in a spherical cap covering, Zong \cite{zong1} showed that 
\begin{equation}\label{eq:blow-up}
\left|\xi\right|\leq \left\lfloor \frac{7d}{\ln \epsilon}\right\rfloor^{\alpha14^{d}d^{2d+3}(\ln\epsilon)^{-d}}, 
\end{equation}
\noindent where $c$ is an absolute constant. Since Zong's construction does not provide much room for improving the above estimate, better constructions are needed to reduce the size of $\xi$, while at the same time keeping $\epsilon$ sufficiently small. 

\begin{problem}\label{zongopen}
Develop a computationally efficient procedure for constructing $\epsilon$-nets of ${\cal{K}}^{d}$ of small cardinality. 
\end{problem}

Addressing the above problem would be a critical first step in implementing Zong's program. Wu \cite{wu-chinese} (also see \cite{he}) has recently proposed two variants of $\gamma_{m}(\cdot)$ that can be used in Zong's program instead. However, the challenges and implementation issues remain the same.

\section*{Acknowledgments}
 The first author is partially supported by a Natural Sciences and Engineering Research Council of Canada Discovery Grant. The second author is supported by a Vanier Canada Graduate Scholarship (NSERC) and Alberta Innovates Technology Futures (AITF).

\bigskip

\normalsize
\noindent K\'aroly Bezdek \\
\small{Department of Mathematics and Statistics, University of Calgary, Canada}\\
\small{Department of Mathematics, University of Pannonia, Veszpr\'em, Hungary\\
\small{\texttt{E-mail:bezdek@math.ucalgary.ca}}

\normalsize

\bigskip
\noindent and
\bigskip

\noindent Muhammad A. Khan \\
 \small{Department of Mathematics and Statistics, University of Calgary, Canada}\\
 \small{\texttt{E-mail:muhammkh@ucalgary.ca}}

\end{document}